\documentclass[11pt]{amsart}
\usepackage{amsxtra}
\usepackage{latexsym}
\usepackage{amsfonts}
\usepackage{amssymb}
\theoremstyle{plain}

\newcommand{\cleqn}{\setcounter{equation}{0}}
\newcommand{\clth}{\setcounter{theorem}{0}}
\newcommand {\sectionnew}[1]{\section{#1}\cleqn\clth}

\newtheorem{theorem}{Theorem}[section]
\newtheorem{lemma}[theorem]{Lemma}
\newtheorem{definition-theorem}[theorem]{Definition-Theorem}
\newtheorem{proposition}[theorem]{Proposition}
\newtheorem{corollary}[theorem]{Corollary}
\newtheorem{definition}[theorem]{Definition}
\newtheorem{example}[theorem]{Example}
\newtheorem{remark}[theorem]{Remark}
\newtheorem{notation}[theorem]{Notation}
\newcommand \bth[1] { \begin{theorem}\label{t#1} }
\newcommand \ble[1] { \begin{lemma}\label{l#1} }

\newcommand \bpr[1] { \begin{proposition}\label{p#1} }
\newcommand \bco[1] { \begin{corollary}\label{c#1} }
\newcommand \bde[1] { \begin{definition}\label{d#1}\rm }
\newcommand \bex[1] { \begin{example}\label{e#1}\rm }
\newcommand \bre[1] { \begin{remark}\label{r#1}\rm }

\newcommand \bnota[1] { \begin{notation}\label{n#1}\rm }

\renewcommand {\eth} { \end{theorem} }
\newcommand {\ele} { \end{lemma} }

\newcommand {\epr} { \end{proposition} }
\newcommand {\eco} { \end{corollary} }
\newcommand {\ede} { \end{definition} }
\newcommand {\eex} { \end{example} }
\newcommand {\ere} { \end{remark} }

\newcommand {\enota} { \end{notation} }
\newcommand \thref[1]{Theorem \ref{t#1}}
\newcommand \leref[1]{Lemma \ref{l#1}}

\newcommand \coref[1]{Corollary \ref{c#1}}
\newcommand \deref[1]{Definition \ref{d#1}}
\newcommand \exref[1]{Example \ref{e#1}}
\newcommand \reref[1]{Remark \ref{r#1}}
\newcommand \lb[1]{\label{#1}}
\def \d {{\partial}}   

\def \Rset {{\mathbb R}}         

\def \O  {{\mathcal{O}}}

\def \V {{\mathcal{V}}}

\def \T {{\mathcal{T}}}
\def \G {{\mathcal{G}}}

\def \ra  {\rangle}
\def \la  {\langle}
\def \lara {\la \, , \, \ra}           
\def \lra {\longrightarrow}

\def \Map {\longmapsto}

\def \wt {\widetilde}

\def \hs {\hspace{.2in}}

\def \Ad { {\mathrm{Ad}} }

\def \g  {\mathfrak{g}}   
\def \h  {\mathfrak{h}}

\def \d  {\mathfrak{d}}

\def \fg  {\mathfrak{g}}

\def \p  {\mathfrak{p}}

\def \z  {\mathfrak{z}}

\def \l {\mathfrak{l}}

\DeclareMathOperator \ad { {\mathrm{ad}} }

\DeclareMathOperator \End { {\mathrm{End}} }

\DeclareMathOperator \Hom { {\mathrm{Hom}} }

\def \lrw {\longrightarrow}
\def \Map {\longmapsto}
\def \beqa {\begin{eqnarray*}}
\def \eeqa {\end{eqnarray*}}

\def \tpi {\tilde{\pi}}
\def \top {{\rm top}}

\def \ue {{\underline{e}}}

\def \piG {\pi_{{\scriptscriptstyle G}}}
\def \piGL {\pi_{{\scriptscriptstyle {G, \Lambda}}}}
\def \piDGL {\pi_{{\scriptscriptstyle {D/G, \Lambda}}}}
\def \piP {\pi_{{\scriptscriptstyle P}}}
\def \ud {\underline{d}}
\def \piDG {\pi_{{\scriptscriptstyle {D/G}}}}

\begin{document}
\title[A note on Poisson homogeneous spaces]
{A note on Poisson homogeneous spaces}
\author[Jiang-Hua Lu]{Jiang-Hua Lu}
\address{
Department of Mathematics   \\
The University of Hong Kong \\
Pokfulam Road               \\
Hong Kong}
\email{jhlu@maths.hku.hk}
\date{}
\begin{abstract}
We identify the cotangent bundle Lie algebroid of a 
Poisson homogeneous space $G/H$ of a Poisson Lie group
$G$ as a quotient of a transformation
Lie algebroid over $G$. As applications, we describe the modular
vector fields of $G/H$, and we identify the
Poisson cohomology of $G/H$ with  coefficients in
powers of its canonical line bundle with
relative Lie algebra cohomology of the Drinfeld
Lie algebra associated to $G/H$. We also construct a 
Poisson groupoid over $(G/H, \pi)$ which is symplectic
near the identity section. 
This note serves as preparation for forthcoming papers, in which
we will compute explicitly the Poisson cohomology and study their symplectic groupoids
for
certain examples of Poisson homogeneous spaces related to 
semi-simple Lie groups. 
\end{abstract}
\maketitle
\sectionnew{Introduction}\lb{intro}

The cotangent bundle of a Poisson manifold $(P, \pi)$ 
is naturally a Lie algebroid \cite{vaisman}  called
the cotangent bundle Lie algebroid of $(P, \pi)$ and denoted by 
$T^*(P, \pi)$. Let
$K_P=\wedge^\top T^*P$ be the canonical line bundle over $P$.
Then the   Lie algebroid $T^*(P, \pi)$ has
a natural representation on $K_P$.
The Poisson cohomology  of $(P, \pi)$ as defined in 
\cite{li:poi}, the Poisson homology of $(P, \pi)$
as defined in
\cite{brylinsky}, and the twisted Poisson cohomology 
of $(P, \pi)$ as defined in \cite{e-lu-we}, 
can be regarded as the Lie algebroid cohomology of $T^*(P, \pi)$
with coefficients in, respectively, the trivial line bundle, $K_P$ and $K_P^2$ 
(see \cite{e-lu-we, vaisman, xuping}). 
In general, one can consider the Lie algebroid cohomology of
$T^*(P, \pi)$ with coefficients in $K_P^N$ for any integer $N$, which we will denote
by $H^\bullet(P, \pi; K_P^N)$ and refer to as 
{\it generalized Poisson cohomology} of $(P, \pi)$.
A symplectic groupoid of $(P, \pi)$ is a Lie groupoid over $P$ with 
Lie algebroid $T^*(P, \pi)$ and a compatible
symplectic structure \cite{we:spoid}.

This note concerns the cotangent bundle Lie algebroids of Poisson homogeneous spaces 
of a Poisson
Lie group $(G, \piG)$. More precisely, by a theorem of Drinfeld
\cite{dr:homog}, each Poisson homogeneous
space $(G/H, \pi)$ of $(G, \piG)$ corresponds to a Lie subalgebra $\l$
of the double Lie algebra $\d$ of $(G, \piG)$. In this note, we identify the
cotangent bundle Lie algebroid of $(G/H, \pi)$ with a {\it quotient} of 
the transformation Lie algebroid $G \rtimes_\lambda \l$ over $G$ associated to
an infinitesimal
action $\lambda$ of $\l$ on $G$. We also identify the representation of
$T^*(G/H, \pi)$ on $K_{G/H}$ with a {\it quotient representation}
of $G \rtimes_\lambda \l$ (see $\S$\ref{quot-trans} for the detail). 

We give two applications. First,  
for any integer $N$, we identify 
the generalized Poisson 
cohomology $H^\bullet\left(G/H, \pi; K_{G/H}^N\right)$ with Lie algebra cohomology of 
$\l$ relative to $H$ with coefficients in $C^\infty(G)_N$, the space
of smooth functions on $G$ together with an $(\l, H)$-module 
structure that depends on $N$
(see \coref{cohom-N} for detail). 
We also discuss the canonical pairing between $H^\bullet\left(G/H, \pi; K_{G/H}^N\right)$
and $H^\bullet\left(G/H, \pi; K_{G/H}^{2-N}\right)$
as a pairing on relative Lie algebra cohomology of $\l$, and we compute the
modular vector fields of $(G/H, \pi)$.
The identifications of the Poisson cohomology and homology (i.e., when $N = 0$ and $N = 1$)
with relative Lie algebra cohomology of $\l$ have been 
established in \cite{lu:homog} and \cite{so} but by
different methods. 

As a second application, we construct a Poisson groupoid $\Gamma$
over $(G/H, \pi)$ that is symplectic near the identity section, and we give
conditions and examples when it is symplectic. The groupoid 
structure on $\Gamma$ is a quotient of a transformation groupoid over $G$
(see Mackenzie's book \cite{mackenzie2} for a general treatment of quotients of groupoids),
while the Poisson structure on $\Gamma$ is obtained by reduction of
a quasi-Poisson manifold  by an action of a quasi-Poisson Lie group,
a theory developed by Alekseev and Kosmann-Schwarzbach in \cite{anton-yvette}.
In the special case when $(G, \piG)$ is complete and when
$H$ is a Poisson Lie subgroup of $(G, \piG)$
with $\pi$ being the projection of $\piG$ to $ G/H$, a symplectic groupoid of
$(G/H, \pi)$ was
constructed by P. Xu  in \cite{xuping-oid}.

There are many examples of Poisson homogeneous spaces associated to
semi-simple Lie groups, and they are in general not
of the type $G/H$ with $H$ being
a Poisson Lie subgroup. See \cite{e-lu-reallag, e-lu-cplxlag, lu-yakimov-3} 
for studies of certain varieties which can serve as moduli spaces of
Poisson homogeneous spaces. 
In forthcoming papers, we will use results from this note to 
compute explicitly the Poisson cohomology and study their symplectic groupoids
for
certain examples of Poisson homogeneous spaces treated in 
\cite{e-lu-reallag, e-lu-cplxlag, lu-yakimov-3}. 
Such examples included flag varieties of complex semi-simple
groups \cite{e-lu-reallag} and semi-simple Riemannian symmetric spaces \cite{foth-lu:symmetric}
(see \exref{ex-semi-simple}).

\subsection{Notation}\lb{notation}
For a smooth manifold $P$, the tangent and cotangent bundles of $P$
are denoted by $TP$ and $T^*P$ respectively. For an integer $0 \leq k \leq \dim P$,
$\V^k(P)$ and $\Omega^k(P)$ will denote respectively the spaces of 
smooth $k$-vector fields and smooth $k$-forms on $P$, and 
\[
\V(P) = \oplus_{k=0}^{\dim P} \V^k(P) \hspace{.2in} {\rm and} 
\hspace{.2in} \Omega(P) = \oplus_{k=0}^{\dim P} \Omega^k(P).
\]
If
$P$ and $Q$ are smooth manifolds and
$F: P \to Q$ is a smooth map, $F_*$ will denote the induced map $TP \to TQ$. 

For a vector bundle $A$ over   $P$, 
$\Gamma(A)$ will denote the space of smooth sections of $A$.
If $V$ is an $n$-dimensional vector space, $\wedge^\top V$ always denotes 
$\wedge^n V$. Let $V^*$ be the dual space of $V$.
For $x \in \wedge^kV$ and $\xi \in \wedge^jV^*$ with 
$k \leq j$, $\iota_x\xi \in \wedge^{j-k}V^*$ is defined by
$(\iota_x \xi, y) = (\xi, x \wedge y)$ for all $y \in \wedge^{j-k}V$.
Unless otherwise specified, all vector spaces are real.

For a Lie group $G$ and $g \in G$, $l_g$ and $r_g$ denote respectively the 
left and right translation on $G$ by $g$. The identity element of a group is always denoted
by $e$.

\subsection{Acknowledgement} We thank K. Mackenzie for references on quotients of 
Lie algebroids and groupoids and Bing-Kwan So for 
helpful discussions. Research for this paper was partially supported by 
HKRGC grants 701603, 703304, and 
the HKU Seed Funding for basic research. 

\sectionnew{Some basic facts on Lie algebroids}
\lb{sec-lieoid}

We
refer to \cite{mackenzie1, mackenzie2} for details
on the facts reviewed in this section.

\subsection{Lie algebroids and Lie algebroid cohomology}
\lb{lieoid-cohom}
Recall that a Lie algebroid over a manifold $P$ is a vector bundle $A$ over $P$
together with a vector bundle homomorphism $\rho_A: A \to TP$ and a Lie bracket $[ \, , \, ]$
on $\Gamma(A)$ such that 

1) $[fa_1, a_2] = f[a_1, a_2] -\rho_A(a_2)(f) a_1$ for all $f \in C^\infty(P)$ and $a_1, a_2 \in
\Gamma(A)$;

2) $\rho_A[a_1, a_2] = [\rho_A(a_1), \rho_A(a_2)]$ for all $a_1, a_2 \in
\Gamma(A)$.

\noindent
Let $A$ be a Lie algebroid over $P$. A {\it representation} of  $A$   
on a vector bundle $E$ over $P$ is
an $\Rset$-bilinear map
$D:  \Gamma(A) \times \Gamma(E) \rightarrow \Gamma(E): 
(a, s) \mapsto D_a s$, 
such that for any $a, b \in \Gamma(A)$, $s \in \Gamma(E)$, and $f \in C^\infty(P)$,
\begin{eqnarray*}
&1) & D_{fa} s = f D_{a} s; \\
&2) & D_a(fs) = fD_a s + (\rho(a)f)s;\\
&3) & D_a(D_b s) - D_b(D_a s)= D_{[a, b]}s.
\end{eqnarray*}
The {\it trivial representation} of $A$ is the one on the
trivial line bundle $E = P \times \Rset$
given by
$D_a f = \rho(a) (f)$ for $a \in
\Gamma(A)$ and $f \in \Gamma(E) \cong C^\infty(P)$.
One has the natural notion of tensor products and duals of 
representations of $A$. In particular, a representation $D$ of $A$ on a 
line bundle $L$ gives rise to  a representation of
$A$ on the $N$-th power $L^N$ of $L$ for any integer $N \geq 0$. For a negative integer
$N$, we use the natural identification between $L^N$ and$(L^{-N})^*$ and thus have a 
representation of $A$ on
$L^N$ as well. 

For a representation $D$ of $A$ on $E$, and for $k \geq 0$, define 
\begin{eqnarray*}
d_{A, E}: \; \Gamma(\Hom(\wedge^k A, E) ) &\lra& \Gamma(\Hom(\wedge^{k+1}A, E))\\
(d_{A, E}\phi)(a_1, a_2, \cdots, a_{k+1}) &= &\sum_{j=1}^{k+1} 
(-1)^{j+1}D_{a_j} \phi(a_1, \cdots, \hat{a}_j,
\cdots, a_{k+1}) \\ 
& & + \sum_{i<j} (-1)^{i+j} \phi([a_i, a_j], \cdots, \hat{a}_i, \cdots 
\hat{a}_j, \cdots, a_{k+1})
\end{eqnarray*}
for $a_1, \ldots, a_{k+1} \in \Gamma(A)$.
Then $d_{A, E}^{2} = 0$. 
The cohomology of the
cochain complex  
\[
(\Gamma(\Hom(\wedge A, E)), \;d_{A, E}),
\]
which will be denoted by $H^{\bullet}_{{\rm Lie}}(A; E)$, 
is called the {\it Lie algebroid cohomology of $A$ with coefficients in $E$}. 
When $E$ is the trivial representation, we  set
$H^\bullet(A; E) = H^{\bullet}_{{\rm Lie}}(A)$.

\subsection{Relative Lie algebra cohomology}\lb{relative}
Our reference for this section is \cite{borel-wallach}.
A Lie algebra $\l$ can be regarded  as a Lie algebroid over
a one point space, so for every $\l$-module $V$, we have the
coboundary operators
\[
d_{\l, V}: \; \Hom(\wedge^k \l, V) \lra \Hom(\wedge^{k+1} \l, V),
\hs k \geq 0.
\]
Let
$\h \subset \l$ be a Lie subalgebra,
$H$ a Lie group with Lie algebra $\h$, and $H \to {\rm Aut}(\l): \; h \mapsto \Ad_h$
a group homomorphism integrating the adjoint action of $\h$ on $\l$.

\bde{lHmodules}
An {\it $(\l, H)$-module} is
a topological vector space $V$ which is both an $\l$-module and 
an $H$-module such that

1) for every $v \in V$,
the map $H \to V: h \mapsto hv$ is smooth, and that the
restriction to $\h$ of the action of $\l$ on $V$ coincides with the
one induced from the $H$-action;

2) for every $v \in V, x \in \l$, and $h \in H$, 
$h(x(h^{-1}(v))) = (\Ad_{h}x)(v)$.
\ede

Let $V$ be an $(\l, H)$-module. For $k \geq 0$, let
\[
C^{k}_{\l, H; V} = \left(\wedge^k(\l/\h)^* \otimes V\right)^H,
\]
where the superscript $H$ denotes the subspace of $H$-invariants. 
Identify $(\l/\h)^* \cong \{\xi \in \l^* \mid \xi |_\h = 0\} \subset \l^*$ and 
regard $C^{k}_{\l, H; V}$ as in 
$\wedge^k \l^* \otimes V \cong \Hom(\wedge^k \l, V)$. 
Then
\[
\bigoplus_{k \geq 0} C^{k}_{\l, H; V} \; \subset 
\bigoplus_{k\geq 0} \Hom (\wedge^k \l, V)
\]
is invariant under $d_{\l, V}$. The cohomology of the
cochain complex $(C^{\bullet}_{\l, H; V}, \, d_{\l, V})$,
which will be denoted by
 $H_{{\rm Lie}}^{\bullet}(\l, H; V)$, 
is called the {\it Lie algebra cohomology of
$\l$ relative to $H$ with coefficients in $V$}.

Suppose that $U$ and $V$ are two $(\l, H)$-modules. Then $U \otimes V$
is naturally an $(\l, H)$-module. For any $0 \leq
j, k \leq n= \dim (\l/\h)$, define
\[
C_{\l, H; U}^j \times C_{\l, H; V}^k \lra C_{\l, H; \, U \otimes V}^{j+k}: 
\; \; 
(c_1, c_2) \longmapsto c_1 \otimes c_2 := \phi \wedge \psi \otimes u \otimes v,
\]
where $c_1 = \phi \otimes u, c_2 = \psi \otimes v$ with
$\phi \in \wedge^j(\l/\h)^*$, 
$\psi \in \wedge^{k} (\l/\h)^*$, $u \in U$, and $v \in V$.
It is easy to check that
\begin{equation}\label{dcc}
d_{\l, \, U \otimes V}(c_1 \otimes c_2) = d_{\l, U} (c_1) \otimes c_2 + 
(-1)^{j} c_1 \otimes d_{\l, V}(c_2)
\end{equation}
if $c_1 \in C_{\l, H; U}^j$.  
Assume that
$ \nu \in \left(C_{\l, H; \, U \otimes V}^{n} \right)^*$ is
such that 
\begin{equation}\lb{F-zero}
\nu \left(d_{\l, \, U \otimes V}(C_{\l, \,H;\, U \otimes V}^{n-1})\right) = 0.
\end{equation}
For $0 \leq k \leq n$, define the pairing $( \, , \, )_\nu$ between 
$C_{\l, H; \;U}^{k}$ and $C_{\l, H; \;V}^{n-k}$ by
\[
(c_1, \; c_2)_\nu = \nu(c_1 \otimes c_2).
\]
It follows from \eqref{dcc} that
\[
(d_{\l, U} (c_1), \; c_2)_\nu + (-1)^{k-1} (c_1, \; d_{\l, V}(c_2))_\nu = 0
\]
for all $c_1 \in C_{\l, H; U}^{k-1}$ and $c_2 \in C_{\l, H; V}^{n-k}$.
Thus $(\, , \, )_\nu$ induces a well-defined pairing, still denoted by
$( \, , \, )_\nu$,  between 
$H^{k}_{{\rm Lie}}(\l, H; U)$ and $H^{n-k}_{{\rm Lie}}(\l, H; V)$ 
for every $0 \leq 
k \leq n$.

\subsection{Quotients of transformation Lie algebroids}\lb{quot-trans}

Let again 
$\l$ be a Lie algebra, $\h \subset \l$  a Lie subalgebra,
$H$  a Lie group with Lie algebra $\h$, and $H \to {\rm Aut}(\l): h \to \Ad_h$
a group homomorphism
integrating the adjoint action of $\h$ on $\l$.

\bde{lHspace}
An {\it $(\l, H)$-space} is a smooth manifold $M$ 
together with a Lie algebra homomorphism  
$\lambda: \l \to \V^1(M)$ and a right action of $H$ on $M$
such that 

1) the restriction of $\lambda$ on $\h$ coincides with the infinitesimal
action of $\h$ on $M$ induced by the right $H$-action, and 

2) for all $m \in M, x \in \l$ and $h \in H$,
$\lambda_{x}(mh) = h_* \lambda_{\Ad_h x}(m)$, where $h_*$ is the
differential of the map $h: M \to M: m_1 \mapsto m_1h$ for $m_1 \in M$.

We will sometimes denote an $(\l, H)$-space by the pair 
$(M, \lambda)$ without explicitly mentioning the action of $H$ on $M$.
\ede

Let  $(M, \lambda)$ be an $(\l, H)$-space. Using
the action $\lambda$ of $\l$ on $M$, one can
form the 
{\it transformation Lie
algebroid} $M \rtimes_\lambda \l$ over $M$, which is the trivial vector 
bundle $M \times \l$ over $M$ with the anchor map
\[
M \times \l \lra TM: \; (m, x) \longmapsto \lambda_x(m), \;
m \in M, x \in \l,
\]
and the Lie bracket $[ \, , \, ]_{M \rtimes_\lambda \l}$
on $\Gamma(M \times \l) \cong C^\infty(M, \l)$ determined by
\[
[\bar{x}_1, \bar{x}_2]_{M \rtimes_\lambda \l} = \overline{[x_1, x_2]},
\]
where for $x \in \l$, $\bar{x}$ is the constant function on $M$ with value $x$.

Assume in addition that the $H$-action on $M$ is free and proper so
that the quotient $M/H$ is a smooth manifold.
Consider 
the associated vector bundle $A = M \times_H (\l/\h)$
over  $M/H$, where $h \in H$ acts on 
$\l/\h$ by $\Ad_h$. Points in $A$ will be denoted by
$[m, x + \h]$, where $m \in M$ and $x \in \l$.
Note that
\[
\rho_A: \; A \lra T(M/H): \; \; [m, \; x+\h] \Map q_*\lambda_x(m)
\]
is a well-defined bundle map,
where $q: M \to  M/H$ is the natural projection, and 
\begin{eqnarray*}
\Gamma(A) &=& C^\infty(M, \l/\h)^H \\
& = &  \{a \in C^\infty(M, \l/\h) \mid
a(mh) = \Ad_{h^{-1}} a(m), \; 
\forall m \in M, h \in H\}.
\end{eqnarray*}
Let 
\[
\Gamma(M \rtimes_\lambda \l)^H = \{a \in C^\infty(M, \l) \mid
a(mh) = \Ad_{h^{-1}} a(m), \; 
\forall m \in M, h \in H\}.
\]
For 
$a_1, a_2 \in \Gamma(A)$, let $\tilde{a}_1, \tilde{a}_2 \in 
\Gamma(M \rtimes_\lambda \l)^H $ be such that 
$\p (\tilde{a}_1) = a_1$ and $\p (\tilde{a}_2) = a_2$, where
$\p: M \rtimes_\lambda \l \to A: (m, x) \mapsto [m, x + \h]$ 
is the natural vector bundle projection.  
Define $[a_1, a_2] \in \Gamma(A)$ by
\begin{equation}
\lb{bra-quotient}
[a_1, a_2] = \p([\tilde{a}_1, \tilde{a}_2]_{M \rtimes_\lambda \l}).
\end{equation}
The proof of the following lemma is omitted since it is straightforward.

\ble{lieoid-quot}
Formula \eqref{bra-quotient} is a well-defined Lie bracket
on $\Gamma(A)$. With the Lie bracket in
\eqref{bra-quotient} on $\Gamma(A)$ and $\rho_A$ as the anchor map,
$A$ is a 
Lie algebroid over $M/H$. Moreover, the bundle map
$\p: M \rtimes_\lambda \l \to A$ is a Lie algebroid
morphism.
\ele

\bde{quotient}
The Lie algebroid $A$ in \leref{lieoid-quot} is called 
the {\it $H$-quotient} of the transformation Lie algebroid $M \rtimes_\lambda \l$
and will be denoted by $M \rtimes_{\lambda, H} (\l/\h)$.
\ede

\bex{ex-GH}
If $G$ is a Lie group and $H \subset G$ a closed subgroup, the tangent
bundle Lie algebroid $T(G/H)$ is a quotient by $H$ of the tangent bundle
Lie algebroid $TG$.
A more general discussion on quotients of Lie algebroid can be
found in \cite[Chap. 4]{mackenzie2}.
\eex

We now turn to a special class of representations of $M \rtimes_{\lambda, H} (\l/\h)$ that arise 
from representations of $M \rtimes_\lambda \l$.

\bde{lHbundles}
An {\it $(\l, H)$-vector bundle} is an $H$-equivariant
 vector bundle $E$ over
an $(\l, H)$-space $(M, \lambda)$ together with a representation of $\l$
on $\Gamma(E)$ such that

1) $x \cdot (fs) = \lambda_x(f) s + f (x \cdot s)$ for all $x \in \l,
f \in C^\infty(M)$, and $s \in \Gamma(E)$;

2) the $\l$-action and the $H$-action on $\Gamma(E)$ induced from
the $H$-action on $E$ make $\Gamma(E)$ into an $(\l, H)$-module
(see \deref{lHmodules}).
\ede

Let $E$ be an $(\l, H)$-vector bundle over $M$ such that the $H$-action on $M$
is free and proper. 
One then has the representation $\wt{D}$ of $M \rtimes_\lambda \l$ on $E$ given by
\[
(\wt{D}_{b} s)(m) = (b(m) \cdot s)(m), \hs b \in 
\Gamma(M \rtimes_\lambda \l) = C^\infty(M, \l), \;
m \in M, s \in \Gamma(E).
\]
Let $E/H$ be the quotient bundle 
over $M/H$ with $\Gamma(E/H) = \Gamma(E)^H$, 
the space of $H$-invariant smooth sections of $E$. 
For $a \in \Gamma(A)$, let $\tilde{a} \in \Gamma(M \rtimes_\lambda \l)^H$ be such that
$\p(\tilde{a}) = a$. It is easy to see that $\wt{D}_{\tilde{a}} s \in \Gamma(E)$
is $H$-invariant for any $s \in \Gamma(E/H) \cong
\Gamma(E)^H$, so we can regard $\wt{D}_{\tilde{a}}s$ as in $\Gamma(E/H)$. Define
\begin{equation}\lb{AEH}
D_a s = \wt{D}_{\tilde{a}} s, \hs s \in \Gamma(E/H) \cong \Gamma(E)^H.
\end{equation}

The proof of the following \leref{AEH2} is straightforward.

\ble{AEH2}
Formula \eqref{AEH} is a well-defined representation of the quotient Lie algebroid
$M \rtimes_{\lambda,H} (\l/\h)$ on $E/H$, and we call it the 
$H$-quotient 
of the representation $\wt{D}$ of $M \rtimes_\lambda  \l$ on $E$.
\ele

\ble{AEH-cohom}
The Lie algebroid cohomology of $A=M \rtimes_{\lambda,H} (\l/\h)$ 
with coefficient in $E/H$ is
isomorphic to the Lie algebra cohomology of $\l$ relative
to $H$ with coefficients in $\Gamma(E)$, i.e.,
\[
H_{{\rm Lie}}^{k}(A; E/H) \; \cong \; H_{{\rm Lie}}^{k}(\l, H; 
\Gamma(E)), \hs \forall k \geq 0.
\]
\ele

\begin{proof} Let $\T$ be the trivial vector bundle over $M$ with
fiber $\l/\h$. Then for every $k \geq 0$, the vector bundle
$\Hom(\wedge^k A, E/H)$ over $ M/H$ is the quotient by $H$ of the
$H$-equivariant vector bundle $\Hom(\wedge^k\T, E)$, so
\begin{equation}\lb{dd}
\Gamma( \Hom(\wedge^k A, E/H)) 
\cong \left(\wedge^k(\l/\h)^* \otimes \Gamma(E)\right)^H
\cong C_{\l, H; \Gamma(E)}^{k}.
\end{equation}
By following the definitions of the Lie algebroid structure on $A$
and the representation of $A$ on $E/H$, it is straightforward to
check that the identifications in \eqref{dd} give an isomorphism
of cochains 
\[
\left( \bigoplus_{k \geq 0} \Gamma( \Hom(\wedge^k A, E/H)), \; d_{A, E/H}
\right) \lra \left(\bigoplus_{k \geq 0} C^{k}_{\l, H; \Gamma(E)}, \; d_{\l, H}
\right).
\]
\end{proof}

\bre{sqroot}
Suppose that $F$ is an $(\l, H)$-line bundle over an $(\l, H)$-space $(M, \lambda)$
and that $E$ is an $H$-equivariant square root of $F$, i.e., $E^2 \cong F$.
Then $E$ is naturally an $(\l, H)$-line bundle with the 
$\l$-action on $\Gamma(E)$ uniquely defined as follows: if $t$ is a nowhere
vanishing local section of $E$, then $x \cdot t = \frac{1}{2}
\frac{x \cdot t^2}{t}$ for any $x \in \l$  (see  \cite{e-lu-we}). Consequently, one has the
quotient representation of the quotient Lie algebroid $A =M \rtimes_{\lambda,H} (\l/\h)$
on $E/H$.
\ere

\sectionnew{Poisson cohomology and modular vector fields}\lb{poi-cohom}
 
\subsection{The cotangent bundle Lie algebroid and Poisson cohomology}
\lb{cotangent-oid}

The cotangent bundle Lie algebroid of a Poisson manifold $(P, \pi)$,
denoted by $T^*(P, \pi)$, is
the cotangent bundle $T^*P$ of $P$ with the anchor map
\[
\tilde{\pi}:\; \; \; T^*P \lra TP: \; \; \; 
\tpi(\alpha)(\beta) = \pi(\alpha, \beta), \hs \alpha, \beta \in \Omega^1(P),
\]
and the Lie bracket $\{ \, ,\, \}_\pi$
on  $\Omega^1(P)$ given by
\begin{equation}\lb{bra}
 \{\alpha, \beta\}_\pi = d (\pi(\alpha, \beta)) + \iota_{\tilde{\pi}(\alpha)}
d\beta - \iota_{\tilde{\pi}(\beta)} d\alpha,
\hs \alpha, \beta \in \Omega^1(P).
\end{equation}
Let
\[
K_P = \wedge^\top T^*P
\]
be the canonical line bundle over $P$.
It is shown in \cite{e-lu-we, xuping} that there
is a representation of the Lie algebroid $T^*(P, \pi)$ on 
$K_P$ given by 
\begin{equation}
\lb{Damu}
D_\alpha \mu = L_{\tpi(\alpha)} \mu + (\pi, d\alpha) \mu = \{\alpha, \mu\}_\pi
-(\pi, d\alpha) \mu = \alpha \wedge d(i_\pi \mu), \; \; \mu \in  \Omega^\top (P),
\end{equation}
where  $\{\, , \, \}_\pi$ is  the Schouten bracket 
on the space $\Omega(P)$ induced from
the bracket in \eqref{bra} on $\Omega^1(P)$.

\bde{canonical-rep}
The representation of $T^*(P, \pi) $ on $K_P$ is called
the {\it canonical representation} of $T^*(P, \pi)$ on $K_P$.
\ede

For any integer $N$, let $K_{P}^{N}$ be the $N$-th power of
$K_{P}$, equipped with the natural extension of the representation
of $T^*(P, \pi)$. When $N$ is negative, we will understand
$K_{P}^{N}$ as $(K_P^{-N})^*$.

\bde{poi-cohomology-N}
For a Poisson manifold $(P, \pi)$ and for any integer $N$, we define
the {\it Poisson cohomology of $(P, \pi)$ with coefficients $K_{P}^{N}$}
to be the Lie algebroid cohomology of  $T^*(P, \pi)$
with coefficients in $K_{P}^{N}$, and we  denote it by
$H^\bullet(P, \pi; K_{P}^{N})$. When $N = 0$, we simply write
$H^\bullet(P, \pi; K_{P}^{N})$ as $H^\bullet(P, \pi)$. The totality of
$H^\bullet(P, \pi; K_{P}^{N})$ for all integers $N$ is called the
{\it generalized Poisson cohomology} of $(P, \pi)$.
\ede

\bre{poi-cohom}
The
Poisson cohomology of $(P, \pi)$ defined in
\cite{li:poi} is  $H^\bullet(P, \pi)$.
It is shown in \cite{e-lu-we, xuping} that 
the Poisson homology of $(P, \pi)$ defined in \cite{brylinsky} 
is isomorphic to 
$H^\bullet(P, \pi; K_P)$. 
In \cite{e-lu-we}, the cohomology $H^\bullet(P, \pi; K_{P}^{2})$
is called the twisted Poisson cohomology of  $(P, \pi)$.
\ere
 
\subsection{The canonical pairing on Poisson cohomology}\lb{pairings-poi}

Suppose that $P$ is compact and oriented. For $0 \leq k \leq n=\dim P$ and 
an integer $N$, set
\[
C_{P, N}^{k} = \Gamma(\Hom(\wedge^k T^*P, \; K_{P}^{N}))\, 
\cong \,\Gamma(\wedge^k TP \otimes K_{P}^{N}).
\]
The natural identifications of bundles
\[
\wedge^k TP \otimes \wedge^{n-k}TP \cong \wedge^n TP, \hs
K_{P}^{N} \otimes K_{P}^{2-N} \cong K_{P}^{2}, \hs
\wedge^n TP \otimes K_{P}^2 \cong K_P
\]
give rise to an identification 
\[
J: \; \left(\wedge^kTP \otimes K_{P}^{N}\right) \otimes
\left(\wedge^{n-k} TP \otimes K_{P}^{2-N}\right) \lra K_P
\]
and thus an $\Rset$-bilinear pairing
\[
( c_1, \; c_2) := \int_P J(c_1, c_2), \hs 
c_1 \in C_{P, N}^{k}, \, c_2 \in C_{P, 2-N}^{n-k}.
\]
A proof similar to that of Theorem 5.1 of \cite{e-lu-we} shows that
$( \,, \, )$ induces a well-defined pairing
between $H^{k}(P, \pi; K_{P}^{N})$ and $H^{n-k}(P, \pi; K_{P}^{2-N})$. 
We will refer to $( \, , \, )$ the {\it canonical pairing} on the 
generalized Poisson 
cohomology of $(P, \pi)$.

\subsection{Modular vector fields}\lb{modular}
Let $(P, \pi)$ be an orientable Poisson manifold, and let 
$\mu$ be a volume form of $P$. 
The {\it modular vector field of
$\pi$ with respect to $\mu$} (see \cite{we:modular})
is defined to be the vector field $\theta_\mu$
on $P$ such that
\[
D_\alpha \mu = (\theta_\mu, \alpha) \mu, \; \; \; \forall
\alpha \in \Omega^1(P),
\]
where $D_\alpha \mu \in \Omega^\top(P)$
is given in \eqref{Damu}.
For an integer $N$, set $d_N=d_{T^*P, K_{P}^{N}} \in {\rm End}
(C_{P, N}^{\bullet})$.

\bpr{modular-cohom} Let $N$ be any integer.
For any volume form $\mu$, the action of the modular vector field $\theta_\mu$
on $C_{P, N}^{\bullet} = \oplus_{k \geq 0} C_{P, N}^{k}$ by Lie derivative commutes
with the operator $d_N$. When $N \neq 1$, the induced action of $\theta_\mu$
on $H^{\bullet}(P, \pi; K_{P}^{N})$ is trivial.
\epr

\begin{proof}
Consider the identification
\[
{\mathcal I}: \; \; \; \V^k(P) \lra C_{P, N}^{k}: \; \; V \longmapsto V \otimes \mu^N.
\]
Since $L_{\theta_\mu} \mu = 0$, $L_{\theta_\mu} \circ {\mathcal I} 
= {\mathcal I} \circ L_{\theta_\mu}$.
It is also easy to show (see \cite[Lemma 4.4]{e-lu-we}) that the operator 
$\delta_N:={\mathcal I}^{-1} \circ d_N \circ {\mathcal I}$ is given by
\[
\delta_N: \; \V^k(P) \lra \V^{k+1}(P): \; \; V \longmapsto [\pi, V] + N \theta_\mu 
\wedge V.
\]
Since $L_{\theta_\mu} \pi = 0$, it 
is clear that $L_{\theta_\mu}$ commutes with $\delta_N$. 
Consider the operator 
\[
b_{\mu}: \; \; \V^k(P) \lra \V^{k-1}(P): \, \, \iota_{b_{\mu} V} \mu = (-1)^k d 
(\iota_V \mu).
\]
It is easy to see that $b_{\mu}^{2} = 0$ and that $\theta_\mu = b_\mu \pi$.
Moreover, for $V_1 \in \V^k(P)$ and $V_2 \in \V(P)$,
\begin{eqnarray*}
b_\mu (V_1 \wedge V_2) & = & b_{\mu}(V_1) \wedge V_2 + (-1)^k V_1 \wedge
b_\mu (V_2) + (-1)^k [V_1, V_2]\\
b_\mu [V_1, V_2] & = & [b_{\mu}(V_1), V_2]  + (-1)^{k-1}[ V_1,
b_\mu (V_2)]
\end{eqnarray*}
It follows that $b_\mu \delta_N + \delta_N b_\mu = (1-N)L_{\theta_\mu}$. Thus 
$\theta_\mu$ acts trivially on $H^{\bullet}(P, \pi; K_{P}^{N})$ when 
$N \neq 1$.
\end{proof}

\sectionnew{The cotangent bundle Lie algebroids and generalized
Poisson cohomology of
Poisson homogeneous spaces}
\label{sec_poi-lie}

\subsection{Review on Poisson Lie groups}\lb{subsec-poi-lie}
Recall that \cite{dr:bigbra, soibelman, lu-we:poi,sts:dressing} 
a {\it Poisson Lie group} 
is a Lie group $G$ with a Poisson structure $\piG$ such that
the group multiplication
map $(G, \piG) \times (G, \piG) \to (G, \piG): \; (g, h) \mapsto gh$
is Poisson. Let $(G, \piG)$ be a Poisson Lie group. Then $\piG$ necessarily vanishes
at the identity element $e$ of $G$. Let $\delta: \g \to \wedge^2 \g$ be the
linearization of $\piG$ at $e$. Then the dual map
\[
\delta^*: \; \wedge^2 \g^* \lra \g^*: \; \; \xi \wedge \eta \longmapsto [\xi, \eta]
\]
of $\delta$ defines a Lie bracket on $\g^*$, and the pair
$(\g, \delta)$ becomes a {\it Lie bialgebra} \cite{dr:bigbra}. For $x  \in \g$
and $\xi \in \g^*$, define $\ad_x^* \xi \in \g^*$ and $\ad_\xi^* x \in \g$
by
\[
(ad_x^* \xi, \, y) = (\xi, \, [y, x]) \hspace{.2in} \mbox{and}
\hspace{.2in} (\ad_\xi^* x, \, \eta) = (x, \, [\eta, \xi]), 
\hspace{.2in} \mbox{where}
\;\;y \in \g, \eta \in \g^*.
\]
Let $\d = \g \oplus \g^*$. Then the bracket on $\d$ given by
\[
[x + \xi, \, y + \eta] = [x, y] + \ad_\xi^* y - \ad_\eta^* x +
[\xi, \eta] + \ad_x^* \eta - \ad_y^* \xi,\; \; \; x, y \in \g, \xi, \eta \in \g^*,
\]
is a Lie bracket, and the bilinear form $\lara$ on $\d$
given by
\[
\la x + \xi, \;y + \eta \ra  = (x, \eta) + (y, \xi), \hs x, y \in \g,\; \xi, \eta \in \g^*,
\]
is ad-invariant with respect to $[ \, , \, ]$. 
The pair $(\d, \lara)$ is called the double of
the Lie bialgebra $(\g, \delta)$. 
The adjoint action of $\g$ on $\d$ integrates to 
an action of $G$ on $\d$, still denoted by $\Ad_g: \d \to \d$ for $g \in G$, which is
given by \cite{dr:homog}
\begin{equation}
\lb{Gond}
\Ad_g (x + \xi) =: \Ad_{g} x + \iota_{\Ad_{g^{-1}}^{*} \xi} (r_{g^{-1}} \piG(g)) +
\Ad_{g^{-1}}^{*} \xi,
\end{equation}
where $\Ad_g: \g \to \g$ and $\Ad_{g^{-1}}^{*}:
\g^* \to \g^*$ are  the adjoint and co-adjoint actions
of $g \in G$ on $\g$ and on $\g^*$ respectively. 
A subspace $\l$ of $\d$ is said to be Lagrangian if 
$\la x, y \ra = 0$ for all $x, y \in \l$ and if $\dim 
\l = \dim \g$. 

\subsection{Drinfeld Lagrangian subalgebras}
\lb{subsec-drinfi}
Let $H$ be a closed subgroup of $G$.

\bde{de-poi-homog} \cite{dr:homog} A
$(G, \piG)$-homogeneous Poisson structure on $G/H$ is a bivector field
$\pi$ on $G/H$ such that 1) $\pi$ is Poisson, and 2)
the map
\begin{equation}\lb{eq-sigma-GGH}
\sigma: \; \; (G, \piG) \times (G/H, \pi) \lra (G/H, \pi): \; \; \; 
(g_1, g_2H) \longmapsto
g_1g_2H
\end{equation}
is Poisson. 
\ede

By definition, the map $\sigma$ in \eqref{eq-sigma-GGH} is 
Poisson if and only if
\begin{equation}\lb{pi-well}
\pi(gH) = (\sigma_g)_* \pi(eH) + q_* \piG(g), \hs \forall g \in G,
\end{equation}
where $q: G \to G/H$ is the projection,
and for $g \in G$, $\sigma_g: G/H \to G/H$ is defined by
$ g_1H\to gg_1H$ for $g_1 \in G$. 
Thus, $\pi$ is uniquely determined by $\pi(eH) \in \wedge^2 T_{eH} (G/H)$,
and Conditions 1) and 2) on $\pi$ in \deref{de-poi-homog}
become the following two conditions on 
$\pi(eH) \in \wedge^2 T_{eH} (G/H)$:

(i)  $\pi(eH) = (\sigma_h)_* \pi(eH) + q_* \piG(h)$ for all $h \in H$ (so that
$\pi$ given by \eqref{pi-well} is well-defined); and 

(ii) the bi-vector field $\pi$ on $G/H$ determined by $\pi(eH)$ via \eqref{pi-well} is
Poisson.

\noindent
Let $\h$ be the Lie algebra of $H$. Simple linear algebra arguments
show that there is a one to one correspondence between $\wedge^2(\g/\h)$ 
and the set of Lagrangian subspaces $\l$ of $\d$ such that  $\l \cap \g =\h$.
The explicit correspondence is given by
\begin{equation}\lb{r-l}
\wedge^2 (\g/\h) \ni r \mapsto \l_r:=
\{x + \xi \mid x \in \g, \xi \in \g^*, \xi|_{\h} = 0,
\iota_\xi r = x + \h\}.
\end{equation} 
Identify 
$T_{eH}(G/H) \cong \g/\h$. Then an element
$\pi(eH) \in \wedge^2 T_{eH} (G/H)\cong \wedge^2 (\g/\h)$ corresponds to
the Lagrangian subspace $\l_{\pi(eH)}$ of $\d$.
Drinfeld showed \cite{dr:homog} that Conditions (i) and (ii) on 
$\pi(eH) \in \wedge^2 T_{eH} (G/H)$ 
are respectively equivalent to

(a) $\Ad_h \l_{\pi(eH)} = \l_{\pi(eH)}$ for all $h \in H$, where
$\Ad_h: \d \to \d$ is given in \eqref{Gond}, and 

(b) $\l_{\pi(eH)}$ is a Lie subalgebra of $\d$.

\bde{de-drinfi}
When $(G/H, \pi)$ is a Poisson homogeneous space of $(G, \piG)$,
the Lie subalgebra $\l_{\pi(eH)}$ of $\d$ is called the Drinfeld Lagrangian
subalgebra  associated to $\pi(eH)$.
\ede

Let $(G/H, \pi)$ be a Poisson homogeneous space of $(G, \piG)$. Let
$q$ also denote the   projection  $\g \to \g/\h$. Let $\Lambda \in \wedge^2 \g$
be any element such that 
\begin{equation}\lb{eq-q-Lambda}
q(\Lambda) = \pi(eH) \in \wedge^2 T_{eH} (G/H) \cong
\wedge^2 (\g/\h). 
\end{equation}
The following \leref{le-Lambda} is straightforward to prove \cite{dina}.

\ble{le-Lambda} Conditions (i) and (ii) on
$\pi(eH)$ are equivalent to

1) $\Ad_h \Lambda -\Lambda +(r_{h^{-1}})_* \piG(h) \in \h \wedge \g$ for all $h \in H$;

2) $[\Lambda, \Lambda] + 2\delta(\Lambda) \in \h \wedge \g \wedge \g,$

\noindent
where $[ \, , \, ]$ is the
Schouten bracket on $\wedge \g$ and $\delta: \g \to \wedge^2 \g$ is the linearization of
$\piG$ at $e$ as well as its extension $\delta: \wedge^2 \g \to \wedge^{3} \g$ given by
\[
\delta(x \wedge y \wedge z)= \delta(x) \wedge y \wedge z -x \wedge \delta(y) \wedge z + 
x \wedge y \wedge \delta(z), \hs x, y, z \in \g.
\]
\ele

For $\Lambda \in \wedge^2 \g$ as in \eqref{eq-q-Lambda},
define the bi-vector field $\pi_\Lambda$ on $G$ by
\begin{equation}\lb{eq-piLambda}
\pi_\Lambda = \Lambda^l + \piG,
\end{equation}
where $\Lambda^l$ is the left invariant bi-vector field on $G$ with value $\Lambda$ at
$e$. Condition 1) on $\Lambda$ in \leref{le-Lambda}
implies that $q_* \pi_\Lambda$ is a
well-defined bi-vector field on $G/H$. In fact,
\[
q_* \pi_\Lambda = \pi.
\]
Let $\Lambda \xi = \iota_\xi \Lambda$ for $ \xi \in \g^*.$
The Drinfeld Lagrangian subalgebra $\l_{\pi(eH)}$ is  also given by 
\begin{equation}\lb{eq-l-Lambda}
\l_{\pi(eH)} =\{x + \Lambda \xi + \xi \mid x \in \h, \xi \in \g^*,
\xi|_\h = 0\}.
\end{equation}
 
\bre{Lambda}
Although the bi-vector field $\pi_\Lambda$ on $G$
is not necessarily Poisson, we can
still define the skew-symmetric bracket $\{ \, ,\, \}_{\pi_\Lambda}$
on $\Omega^1(G)$ by replacing $\pi$ by $\pi_\Lambda$ in
 \eqref{bra}. Moreover, 
the space of left invariant $1$-forms on $G$ is invariant under 
$\{ \, ,\, \}_{\pi_\Lambda}$. In fact, it is easy to show that
\[
\{\xi^l, \eta^l\}_{\pi_\Lambda} = ([\xi, \eta]_\Lambda)^l,
\hs \xi, \eta \in \l,
\]
where for $\zeta \in \g^*$, $\zeta^l$ is the left invariant $1$-form on 
$G$ with value $\zeta$ at $e$, and 
\begin{equation}\lb{bra-Lambda}
[\xi, \eta]_\Lambda \; \stackrel{{\rm def}}{=} \;
[\xi, \eta] + \ad_{\Lambda \xi}^* \eta - \ad_{\Lambda \eta}^* \xi, \hs 
 \xi, \eta \in \fg^*.
\end{equation}
\ere

\ble{le-h0}
Let $\h^0 = \{\xi \in \g^* \mid \xi|_\h = 0\}$. Then 
$[\xi, \eta]_\Lambda \in \h^0$ for all $\xi, \eta \in \h^0$.
\ele

\begin{proof}
The condition $\Ad_h \l_{\pi(eH)} = \l_{\pi(eH)}$
for all $h \in H$ implies that $[x, \l_{\pi(eH)}] \subset \l_{\pi(eH)}$
for all $x \in \h$, so
$[x, \Lambda \xi + \xi] \in \l_{\pi(eH)}$ for all $x \in \h$ and $\xi \in \h^0$,
from which it follows that $[\xi, \eta]_\Lambda \in \h^0$ for all $\xi, \eta \in \h^0$.
See also \cite{dina}.
\end{proof}

Let $\chi_{\h^0, \Lambda} \in (\h^0)^*$ be defined by 
\[
\chi_{\h^0, \Lambda}(\xi) = {\rm tr}(T_\xi ), \hs \xi 
\in \h^0,
\]
where $T_\xi \in \End (\h^0):  \eta \mapsto [\xi, \eta]_\Lambda$ for $\xi, \eta \in \h^0$.
Let $\chi_\l \in \l^*,  \chi_\g \in \g^*$, and $\chi_{\g^*} \in \g$ 
be the adjoint characters of $\l, \g$ and $\g^*$ respectively. 
Let $b\Lambda = \sum_i [x_i, y_i] \in \g$ if $\Lambda = \sum_i
x_i \wedge y_i$. We now prove a fact that will be used in 
the proof of \thref{cotangent-poi-homog}.

\ble{Lambda-chi} 
For every $\xi \in \h^0$, 
\begin{equation}\lb{chi-xi-0}
\chi_{\h^0, \Lambda}(\xi) + (b \Lambda, \; \xi) = 
\frac{1}{2}\left(\chi_\l(\Lambda \xi + \xi)
-\chi_\g (\Lambda \xi) + \chi_{\g^*}(\xi)\right).
\end{equation}
\ele

\begin{proof} For $\xi \in \g^*$, consider the operator 
$T_\xi \in \End(\g^*):  T_\xi (\eta) = [\xi, \; \eta]_\Lambda,$
and define $\chi_{\g^*, \Lambda}(\xi) =
{\rm tr}(T_\xi \in \End(\g^*))$. It is easy to see that
\[
\chi_{\g^*, \Lambda}(\xi) 
= \chi_{\g^*}(\xi) -\chi_\g(\Lambda \xi) - 2(b\Lambda,
\xi), \hspace{.2in} \xi \in \g^*.
\]
For $\xi \in \h^0$, since $T_\xi(\h^0) \subset \h^0$, we have an induced map
$T_\xi \in \End(\g^*/\h^0)$.
Define $\chi(\xi) = {\rm tr}(T_\xi 
\in \End(\g^*/\h^0))$ for $\xi \in \h^0$. Then  
\begin{equation}\lb{chi-xi-1}
\chi_{\h^0, \Lambda}(\xi) = \chi_{\g^*, \Lambda}(\xi) - \chi(\xi)
=\chi_{\g^*}(\xi) -\chi_\g(\Lambda \xi) - 2(b\Lambda, \xi) - \chi(\xi)
\end{equation}
for all $\xi \in \h^0$.
On the other hand, consider the embedding $\kappa: \h^0 \hookrightarrow \l$ by
$\xi \mapsto \Lambda  \xi + \xi$, and let $\p_\h: \l \to \h$ be the
projection with respect to the decomposition $\l = \h + \kappa(\h^0)$.
For $\xi \in \h^0$, let $S_\xi \in \End(\h)$ be the operator
$S_\xi(x) = \p_\h [\Lambda \xi + \xi, \; x]$ for $ x \in \h.$
Then
\[
\chi_{\h^0, \Lambda}(\xi) = \chi_\l(\Lambda \xi + \xi) -{\rm tr}(S_\xi 
\in \End (\h)), \hs \forall \xi \in \h^0.
\]
By identifying $\g^*/\h^0 \cong \h^*$, one can show that 
$-S_{\xi}^{*} = T_\xi \in \End(\g^*/\h^0)$, and so
${\rm tr}(S_\xi \in \End (\h)) = -\chi(\xi)$ for all $\xi \in \h^0$. 
Thus
\begin{equation}\lb{chi-xi-2}
\chi_{\h^0, \Lambda}(\xi) = \chi_\l(\Lambda \xi + \xi) + \chi(\xi),
\hs \forall \xi \in \h^0.
\end{equation}
Adding \eqref{chi-xi-1} and \eqref{chi-xi-2}, we get \eqref{chi-xi-0}.
\end{proof}

\subsection{The cotangent bundle Lie algebroid of $(G/H, \pi)$}
\lb{cotangent-poi-homog}
Let $(G, \piG)$ be a Poisson Lie group.
For $x \in \g$ and $\xi \in \g^*$, let $x^l$ (resp. $\xi^l$)
be the 
left invariant vector field (resp. $1$-form) on $G$
with value $x$ (resp. $\xi$) at $e$.
Then \cite{lu:homog} the map
\begin{equation}
\lb{lambda}
\lambda: \; \d \lrw \V^1(G): \; x + \xi \Map 
\lambda_{x +\xi} :=x^l + \tilde{\pi}_{\scriptscriptstyle G}(\xi^l)
\end{equation}
is a Lie algebra homomorphism from $\d$ to the space $\V^1(G)$ of
vector fields on $G$. Let $\p_\g: \d \to \g$ be the projection along $\g^*$. By \eqref{Gond}, we
also have
\begin{equation}\lb{lambda-1}
\lambda_{x + \xi}(g) = (r_g)_* p_\g \Ad_g(x + \xi), \hs \, g \in G, \, x \in \g, \, \xi \in \g^*.
\end{equation}

Let now $(G/H, \pi)$ be a 
$(G, \piG)$-homogeneous Poisson space, and let 
$\l = \l_{\pi(eH)}$ be the 
Drinfeld Lagrangian subalgebra of $\d$ as in \deref{de-drinfi}.
Then $G$, with the right action of $H$ by right translations
and the infinitesimal action of $\l$ by $\lambda$, becomes
an $(\l, H)$-space in the sense of \deref{lHspace}. Let $G \rtimes_\lambda \l$
be the corresponding transformation Lie algebroid over $G$.

\bth{cotangent-poi-homog}
The cotangent bundle Lie algebroid of $(G/H, \pi)$ is isomorphic
to the $H$-quotient  $A=G \rtimes_{\lambda,H} (\l/\h)$ of the 
transformation Lie algebroid $G \rtimes_\lambda \l$.
\eth

\begin{proof}
Let $\Lambda \in \wedge^2 \g$ be any element with 
$q(\Lambda) = \pi(eH) \in \wedge^2 T_{eH}(G/H) \cong \wedge^2 (\g/\h)$.
Recall that $\h^0 = \{\xi \in \g^* \mid \xi|_\h = 0\}$. The projection
$\l \to \g^*: x + \xi \to \xi$ gives an $H$-equivariant isomorphism 
$\l/\h \to \h^0$ whose inverse is  $\h^0 \to \l/\h: \xi \mapsto 
\xi + \Lambda \xi + \h$.
 
Using left translations by elements in $G$ and the identification
$T_{eH}^{*}(G/H) \cong \h^0$, we have the vector bundle isomorphism
\begin{equation}\lb{I}
I: \; T^*(G/H) \lra 
G \times_H \h^0 \cong G \times_H (\l/\h). 
\end{equation}
It remains to show that $I$ is a Lie algebroid isomorphism. Recall that
$\pi = q_* \pi_\Lambda$, where   $\pi_\Lambda$
is the bi-vector field on $G$ given by
$\pi_\Lambda = \Lambda^l + \piG$.

Let $n = \dim \h^0$, and let
$\xi_1, \xi_2, \ldots, \xi_n$ be a basis of $\h^0$. 
For $\alpha \in \Omega^1(G/H)$, write
\begin{equation}\lb{eq-qalpha}
q^*\alpha = \sum_{j=1}^{n} f_{\alpha, j} \xi_{j}^{l} \;\in\;
\Omega^1(G), \hs {\rm where} \hs f_{\alpha, j} \in C^\infty(G), \;
j = 1, \ldots, n.
\end{equation}
Then $I(\alpha) = \sum_{j=1}^{n}
f_{\alpha, j} \xi_j  \in  C^\infty(G, \h^0)^H
\cong \Gamma(A)$, 
and $b_\alpha = \sum_{j=1}^{n} f_{\alpha, j} (\Lambda \xi_j + \xi_j)
\in C^\infty(G, \l)^H = \Gamma(G \rtimes_\lambda \l)^H$ is an $H$-invariant lifting of
$I(a)$. Using $q_* \pi_\Lambda =\pi$, one has
\begin{eqnarray*}
\tpi(\alpha) & = & q_* \tpi_\Lambda (q^* \alpha) = q_*\left(
\sum_{j=1}^{n} f_{\alpha, j} \tpi_\Lambda(\xi_{j}^{l})\right)  
 =  q_*\left(
\sum_{j=1}^{n} f_{\alpha, j} ( (\Lambda \xi_j)^l + 
\tilde{\pi}_{\scriptscriptstyle G}(\xi_{j}^{l}))
\right) \\
& = & q_*\left(
\sum_{j=1}^{n} f_{\alpha, j}  \lambda_{\Lambda \xi_j + \xi_j} \right)
 =  \rho_A(I(\alpha)).
\end{eqnarray*}
Thus $I$ maps the anchor map of $T^*(G/H)$ to the anchor map $\rho_A$
of $A$. 

It remains to show that $I\{\alpha, \beta\}_\pi = [I(\alpha), I(\beta)]$
for any $\alpha, \beta \in \Omega^1(G/H)$.
Let $\{ \, , \, \}_{\pi_\Lambda}$ be the skew-symmetric bracket 
on $\Omega^1(G)$ defined by replacing $\pi$ by $\pi_\Lambda$ in
\eqref{bra}. Using again the fact 
that $\pi = q_* \pi_\Lambda$, we have
\begin{eqnarray*}
q^*\{\alpha, \beta\}_\pi & = & \{q^*\alpha, q^* \beta\}_{\pi_\Lambda}
=\sum_{j,k} \{f_{\alpha, j} \xi_j^l, \, 
f_{\beta, k} \xi_k^l\}_{\pi_\Lambda} \\ \\
& = & \sum_{j,k} \left(
f_{\alpha, j} \tpi_\Lambda(\xi_{j}^{l}) (f_{\beta, k})
\xi_{k}^{l} \right)- \sum_{j,k} \left(
f_{\beta, k} \tpi_\Lambda(\xi_{k}^{l}) (f_{\alpha, j})
\xi_{j}^{l} \right)\\ & & + 
\sum_{j,k} \left(f_{\alpha, j} f_{\beta, k} \{\xi_{j}^{l}, 
\xi_{k}^{l} \}_{\pi_\Lambda} \right).
\end{eqnarray*}
Thus, by \reref{Lambda},
\begin{eqnarray*}
I\{\alpha, \beta\}_\pi 
& = &\sum_{j,k} \left(
f_{\alpha, j} \lambda_{\Lambda \xi_j + \xi_j} 
(f_{\beta, k})
\xi_{k} \right)- \sum_{j,k} \left(
f_{\beta, k} \lambda_{\Lambda \xi_k + \xi_k}
 (f_{\alpha, j})
\xi_{j} \right)\\
& & + 
\sum_{j,k} \left(f_{\alpha, j} f_{\beta, k} [\xi_{j}, 
\xi_{k}]_\Lambda \right),
\end{eqnarray*} 
where the bracket $[ \,, \, ]_\Lambda$ on $\g^*$ is defined in
\eqref{bra-Lambda}.
On the other hand, using
\[
b_\alpha = \sum_{j=1}^{n} f_{\alpha, j} (\Lambda \xi_j + \xi_j)
\hs {\rm and} \hs
b_\beta = \sum_{k=1}^{n}  f_{\beta, k} (\Lambda \xi_k + \xi_k)
\]
as $H$-invariant liftings of $I(\alpha)$ and $I(\beta)$ to smooth
sections of $G \rtimes_\lambda \l$, one can compute 
$[I(\alpha), I(\beta)]\in \Gamma(A)$
and see
that $I\{\alpha, \beta\} = 
[I(\alpha), I(\beta)]$. 
This completes the proof that $I$ is a Lie algebroid isomorphism.
\end{proof}
 
\subsection{The canonical representation of $T^*(G/H, \pi)$ on $K_{G/H}$}
\lb{rep-TGH}

Let $(G/H, \pi)$ be a Poisson homogeneous space
of $(G, \piG)$.
Let
\[
E = G \times \wedge^\top \h^0
\]
be the trivial $H$-equivariant line bundle over $G$, where
\[
(g,\,  Y) \cdot h  =(gh, \, \Ad_{h}^*Y), \hspace{.2in} 
g \in G, \,Y \in \wedge^\top \h^0.
\]
Then the identification $I: T^*(G/H) \to G \times_H \h^0$  by left translation
induces an identification $I: K_{G/H} \to E/H$. 
In this section, we show that $E$ is naturally an $(\l, H)$-line bundle
and that the canonical
representation of $T^*(G/H)$
on $K_{G/H} \cong E/H$ can be identified with the $H$-quotient of the
representation of $G \rtimes_\lambda \l$ on $E$,
where $\l$ is the Drinfeld Lagrangian subalgebra of
$\d$ associated to $\pi(eH)$, and $\lambda$ is the  infinitesimal
action of $\l$ on $G$ given in
\eqref{lambda} (see \deref{canonical-rep} and
\leref{AEH2}).

Let $\wedge^\top \l$ be the $1$-dimensional $(\l, H)$-module, 
on which $\l$ acts by the adjoint character
$ \chi_\l$ and $h \in H$ acts by $\Ad_h \in {\rm Aut}(\l)$.
The trivial line bundle over $G$ with fiber $\wedge^\top \l$, still
denoted by $\wedge^\top \l$, is then an $(\l, H)$-line bundle.
Regard $\wedge^\top T^*G$  as an $(l, H)$-line bundle, on which $H$ acts by
right translation and $\l$ acts 
by Lie derivatives via $\lambda$. Set
\[
F = \wedge^\top \l \otimes \wedge^\top T^*G.
\]
Then $F$ is an $(\l, H)$-line bundle.
Clearly, left translation in $G$ gives rise to an $H$-equivariant
trivialization
\[
F \stackrel{\cong}{\lra} G \times (\wedge^\top \l \otimes \wedge^\top \g^*),
\]
where $h \in H$ acts on $G \times (\wedge^\top \l \otimes \wedge^\top \g^*)$
by
\[
(g, \, X \otimes \mu) \cdot h = (gh, \,  (\Ad_{h^{-1}} X) \otimes (\Ad_{h}^* \mu)), 
\hs X \in \wedge^\top \l, \, \mu \in
\wedge^\top \g^*.
\]

\ble{le-EF} $\wedge^\top \l \otimes \wedge^\top \g^* \cong (\wedge^\top \h^0)^2$
as $H$-modules, so $E^2 \cong F$ as $H$-equivariant line bundles over $G$.
\ele
 
\begin{proof} For $V \in \{\h, \l, \h^0, \g^*\}$, let
$\chi_{\scriptscriptstyle{H, V}}$ be the character of the $H$-action on 
$\wedge^\top V$ induced from the adjoint and co-adjoint actions.
It is easy to see that
\[
\chi_{\scriptscriptstyle{H,\l}} = \chi_{\scriptscriptstyle{H, \h}}  
\chi_{\scriptscriptstyle{H, \h^0}} \hs \mbox{and} \hs
\chi_{\scriptscriptstyle{H, \g^*}} =  \chi_{\scriptscriptstyle{H, \h}}^{-1} 
\chi_{\scriptscriptstyle{H, \h^0}}.
\]
Thus $\chi_{\scriptscriptstyle{H, \l}} \chi_{\scriptscriptstyle{H, \g^*}} = 
\chi_{\scriptscriptstyle{H, \h^0}}^2.$
\end{proof}

Since $F$ is an $(l, H)$-line bundle, so is $E$ as a square root of $F$ by \reref{sqroot}. 
In the next
\leref{lf}, we determine the $l$-module structure on $\Gamma(E)$.
Recall that $\chi_\l \in \l^*, \chi_\g \in \g^*$ and $\chi_{\g^*}$ are the adjoint
characters of $\l, \g$, and $\g^*$ respectively.

\ble{lf}
Fix $Y_0 \in \wedge^\top \h^0, \, Y_0 \neq 0$, and write
elements in $\Gamma(E) =C^\infty(G,  \wedge^\top \h^0)$ as 
$f Y_0$ for $f \in C^\infty(G)$. Then the $l$-module structure on 
$\Gamma(E)$ is given by
\[
(x + \xi) \cdot (fY_0) = \left( \lambda_{x+\xi} (f) 
+ \frac{1}{2} \left(\chi_l(x+\xi) -\chi_\g(x) + \chi_{\g^*}(\xi) 
-2(\piG, d\xi^l\right) f \right) Y_0
\]
for any $x + \xi \in \l$ and $f \in C^\infty(G)$. 
\ele 

\begin{proof} Fix non-zero elements $X_0 \in \wedge^\top l$ and $\mu_0 
\in \wedge^\top \g^*$, and let $\mu_0^l$ be the left invariant 
volume form on $G$ with $\mu_0^l(e) = \mu_0$. Then $X_0 \otimes \mu_0^l$
is a nowhere vanishing section of $F$. For $x + \xi \in \l$, one has
\begin{eqnarray*}
(x + \xi) \cdot  (X_0 \otimes \mu_0^l)& = & 
 \chi_\l(x+\xi) X_0 \otimes \mu_{0}^{l} 
+ X_0 \otimes L_{\lambda_{x+\xi}} \mu_{0}^{l} \\
& = &  
(\chi_l(x+\xi) -\chi_\g(x)) X_0 \otimes \mu_{0}^{l} + 
 X_0 \otimes  L_{\tilde{\pi}_{\scriptscriptstyle G}(\xi^l)} \mu_{0}^{l}.
\end{eqnarray*}
By \eqref{Damu},  
\[
L_{\tilde{\pi}_{\scriptscriptstyle G}(\xi^l)}\mu_0^l = \{\xi^l, \, \mu_{0}^{l}\}_{\piG} - 
2(\piG, d\xi^l)
\mu_0^l = (\chi_{\g^*}(\xi)-2(\piG, d\xi^l))\mu_0^l,
\]
from which the formula in \leref{lf} follows.
\end{proof}

By $\S$\ref{quot-trans}, the  $(l, H)$-line bundle structure on $E$ gives rise to a 
representation of  
the transformation
Lie algebroid  $G \rtimes_\lambda  \l$ on $E$ and a representation of the 
$H$-quotient Lie algebroid $A = G\rtimes_{\lambda,H} (\l/\h)$ on $E/H$.

\bth{GH-canonical-rep}
Under the identification 
$I: T^*(G/H, \pi) \cong A = G\rtimes_{\lambda,H} (\l/\h)$ 
of Lie algebroids and the identification
$I: K_{G/H} \cong E/H$ of line bundles, 
the canonical
representation of $T^*(G/H, \pi)$
on $K_{G/H}$ becomes the $H$-quotient
representation of $A$ on 
$E/H$.
\eth

\begin{proof}
Denote by $D$ both the canonical representation of $T^*(G/H, \pi)$ on
$K_{G/H}$ and the quotient representation of
$A$ on $E/H$. We need to show that 
\begin{equation}\lb{DD}
D_{I(\alpha)} I(\mu) = I(D_\alpha \mu),
\hs \forall \; \alpha \in 
\Omega^1(G/H), \; \mu \in \Omega^\top (G/H).
\end{equation}
Let $Y_0 = \xi_{1} \wedge \cdots \wedge \xi_{n} \in \wedge^\top \h^0$,
where $\xi_1, \ldots, \xi_n$ is a basis for $\h^0$, and
write 
\[
q^*\alpha = \sum_{j=1}^{n} f_{\alpha, j} \xi_{j}^{l} \;\in\;
\Omega^1(G) \hs {\rm and} \hs 
q^* \mu = \phi \xi_{1}^{l} \wedge \cdots \xi_{n}^{l} \in \Omega^n(G),
\]
where $f_{\alpha, j} \in C^\infty(G)$ for
$j = 1, \ldots, n$, and $\phi \in C^\infty(G)$. Then 
\[
I(\alpha)  = \sum_{j+1}^{n} f_{\alpha, j} \xi_j \in 
C^\infty(G, \h^0)^H \hs \mbox{and} \hs
I(\mu) = \phi Y_0 \in 
\Gamma(E)^H.
\]
Moreover $b_\alpha := \sum_{j=1}^{n} f_{\alpha, j} (\Lambda \xi_j + \xi_j) \in 
\Gamma(G \rtimes_\lambda \l)^H$ is an $H$-invariant  lifting of $I(\alpha)$ to a 
section of $G \rtimes_\lambda \l$. 
Let $\wt{D}$ be the representation of $G \rtimes_\lambda \l$ on $E$. By \leref{lf}, 
\begin{eqnarray}
\lb{Dbf1}
\wt{D}_{b_\alpha} f_\mu & = & \sum_{j=1}^{n} f_{\alpha, j}
\left(\lambda_{\Lambda \xi_j + \xi_j}
(\phi) - (\piG, d \xi_{j}^{l}) \phi\right) Y_0 \\
\nonumber
& & +\frac{1}{2}\sum_{j=1}^{n} f_{\alpha, j} 
\left(\chi_\l(\Lambda \xi_j + \xi_j)
-\chi_\g (\Lambda \xi_j) + \chi_{\g^*}(\xi_j)\right)\phi Y_0.
\end{eqnarray}
On the other hand, let $Y_{0}^{l}$ be the left invariant $n$-form
on $G$ with $Y_{0}^{l}(e) = Y_0$. Then
\begin{eqnarray*}
q^* D_\alpha \mu 
& = & 
q^* \left( \{\alpha,\,  \mu\}_\pi - 
(\pi, \,d\alpha) \mu \right)=
\{q^* \alpha, \; q^* \mu\}_{\pi_\Lambda} - (\pi_\Lambda, \; d q^* \alpha)
q^* \mu \\
& = & 
\sum_{j=1}^{n} \left( 
\{f_{\alpha, j} \xi_{j}^{l}, \; \phi Y_{0}^{l} \}_{\pi_\Lambda} -
(\pi_\Lambda, \; d(f_{\alpha, j} \xi_{j}^{l} ) ) \phi Y_{0}^{l} \right)\\
& = & 
\sum_{j=1}^{n} \left(
f_{\alpha, j} \tpi_\Lambda(\xi_{j}^{l})(\phi) Y_{0}^{l} + 
\{f_{\alpha, j}\xi_{j}^{l}, \; Y_{0}^{l}\}_{\pi_\Lambda} \phi \right) \\
& & - \sum_{j=1}^{n} (\pi_\Lambda, \; 
df_{\alpha, j} \wedge \xi_{j}^{l} + f_{\alpha, j} d \xi_{j}^{l} ) \phi Y_{0}^{l}\\
& = & \sum_{j=1}^{n} 
f_{\alpha, j} \left(\lambda_{\Lambda \xi_j + \xi_j}
(\phi) - (\piG, \; d \xi_{j}^{l}) \phi\right) Y_{0}^{l} 
\, + \, \{f_{\alpha, j}\xi_{j}^{l}, \; Y_{0}^{l}\}_{\pi_\Lambda} \phi\\ 
& & + \sum_{j=1}^{n} \left( \tpi_\Lambda(\xi_{j}^{l})(f_{\alpha, j}) -f_{\alpha, j} 
(\Lambda^l, \; d\xi_{j}^{l}) \right) \phi Y_{0}^{l}. 
\end{eqnarray*}
Using the properties of the Schouten bracket $\{ \, ,\, \}_{\pi_\Lambda}$
on $\Omega(G)$, one has
\[
\sum_{j=1}^{n} \{f_{\alpha, j} \xi_{j}^{l}, \; Y_{0}^{l}\}_{\pi_\Lambda} =
\sum_{j=1}^{n} \left( f_{\alpha, j} \{\xi_{j}^{l}, \; Y_{0}^{l}\}_{\pi_\Lambda}-
\tpi_\Lambda(\xi_{j}^{l})(f_{\alpha, j}) Y_{0}^{l}\right).
\]
Thus  
\begin{eqnarray}
\lb{Dbf2}
q^* D_\alpha \mu &=&\sum_{j=1}^{n} 
f_{\alpha, j} \left(\lambda_{\Lambda \xi_j + \xi_j}
(\phi) - (\piG, \; d \xi_{j}^{l}) \phi \right) Y_{0}^{l} \\
\nonumber
&&+ \sum_{j=1}^{n} 
f_{\alpha, j} \left(
\{\xi_{j}^{l}, \; Y_{0}^{l}\}_{\pi_\Lambda} \phi -
(\Lambda^l, \; d\xi_{j}^{l})\phi  Y_{0}^{l}\right).
\end{eqnarray}
By \leref{Lambda-chi},
\[
\{\xi, \; Y_{0}^{l}\}_{\pi_\Lambda} -(\Lambda^l, \; d\xi^{l})Y_{0}^{l} = 
\frac{1}{2} 
\left(\chi_\l(\Lambda \xi + \xi)
-\chi_\g (\Lambda \xi) + \chi_{\g^*}(\xi)\right) Y_{0}^{l}, \hs \forall \xi \in \h^0.
\]
Comparing with \eqref{Dbf1} and \eqref{Dbf2}, we see that \eqref{DD} holds.
\end{proof}

\subsection{Poisson cohomology of $(G/H, \pi)$}\lb{GH-cohom}
Let the notation be as in $\S$\ref{cotangent-poi-homog}. For any
integer $N$, since $E$ is a trivial line bundle over $G$, 
$\Gamma(E^N) \cong C^\infty(G)$ as vector spaces. The induced $(l, H)$-module structure
on $C^\infty(G)$ is specified as follows.

\bnota{nota-Cinfty-G-N}
For an integer $N$, denote $C^\infty(G)_N$   the vector space
$C^\infty(G)$ with the following $(\l, H)$-module structure: 
for $x + \xi \in \l, h \in H$ and $f \in C^\infty(G)$,
\begin{eqnarray*}
(x + \xi) \cdot_{\scriptscriptstyle{N}} f & = &\lambda_{x+\xi} (f) 
+ \frac{N}{2} \left(\chi_l(x+\xi) -\chi_\g(x) + \chi_{\g^*}(\xi) 
-2(\piG, d\xi^l)\right) f, \\
h \cdot_{\scriptscriptstyle{N}} f 
&=& \left(\chi_{\scriptscriptstyle{H, \h^0}} (h)\right)^N (f \circ r_h),
\end{eqnarray*}
where $\chi_{\scriptscriptstyle{H, \h^0}} (h) = \det(\Ad_{h^{-1}}^*: \h^0 \to \h^0)$
and $r_h$ is the right translation by $h$.
\enota

We can now identify  the Poisson cohomology of $G/H$ with relative Lie algebra
cohomology. \coref{cohom-N}
follows directly from  \leref{AEH-cohom}, \thref{cotangent-poi-homog}, and
\thref{GH-canonical-rep}.

\bco{cohom-N} For any integer $N$, 
\[
H^\bullet\left(G/H, \pi; K_{G/H}^{N}\right) \cong H^{\bullet}_{{\rm Lie}}
(\l, H; C^\infty(G)_N).
\]
where the left hand side is the generalized Poisson cohomology of $(G/H, \pi)$ and
the right hand side is the Lie algebra cohomology of $\l$ relative to $H$ with coefficients in 
$C^\infty(G)_N$.
\eco

The special case of \coref{cohom-N} when $N = 0$ was proved in \cite{lu:homog}.

\subsection{The pairing on the Poisson cohomology}\lb{pairings}

Assume that $G/H$ is compact and orientable with a fixed orientation, so one has
the map
\begin{equation}\lb{int}
\Omega^\top (G/H) \lra \Rset: \; \omega \Map \int_{G/H} \omega.
\end{equation}
 
By $\S$\ref{pairings-poi}, for any integer $N$ and any 
$0 \leq k \leq n = \dim(G/H)$, there is a well-defined 
pairing $( \, ,\, )$ between $H^{k}\left(G/H, \pi; K^{N}_{G/H}\right)$ and 
$H^{n-k}\left(G/H, \pi; K_{G/H}^{2-N}\right)$. In view of \coref{cohom-N}, we now identify
this pairing with a pairing on the corresponding relative
Lie algebra cohomology spaces. Let the notation be as in 
$\S$\ref{rep-TGH}. Then we have the identifications of $H$-modules:
\begin{eqnarray*}
\wedge^\top (\l/\h)^* \otimes \Gamma(E^N) \otimes
\Gamma(E^{2-N}) & \cong & \wedge^\top (\l/\h)^*  \otimes
\Gamma(E^{2}) \\
&\cong &\wedge^\top (\l/\h)^*  \otimes
\Gamma(F) \\
& \cong & \wedge^\top (\l/\h)^* \otimes 
\wedge^\top \l \otimes \Omega^\top(G)\\
& \cong & \wedge^\top (\l/\h)^* \otimes \wedge^\top \l \otimes 
\wedge^\top \g^* \otimes C^\infty(G)\\
& \cong & \wedge^\top (\l/\h)^* \otimes (\wedge^\top \h^0)^2 \otimes C^\infty(G)\\
& \cong & \wedge^\top \h^0 \otimes C^\infty(G),
\end{eqnarray*}
where we used \leref{le-EF} to identify $\wedge^\top \l \otimes \wedge^\top \g^* 
\cong (\wedge^\top \h^0)^2$ and left translation in $G$ to identify
$\Omega^\top (G) \cong \wedge^\top \g^* \otimes C^\infty(G)$.
Thus we have an identification 
\begin{equation}\lb{II}
(\wedge^\top (\l/\h)^* \otimes \Gamma(E^N) \otimes
\Gamma(E^{2-N}))^H \cong (\wedge^\top (\l/\h) 
\otimes C^\infty(G))^H \cong \Omega^\top (G/H).
\end{equation}
Let $\nu: (\wedge^\top (\l/\h)^* \otimes \Gamma(E^N) \otimes
\Gamma(E^{2-N}))^H  \rightarrow \Rset$ 
be the composition of the  identification in \eqref{II} with the integration map
in \eqref{int}. One checks directly that \eqref{F-zero} holds
and that, under the identifications in \coref{cohom-N}, the canonical pairing 
between  $H^{k}\left(G/H, \pi; K_{G/H}^{N}\right)$ and  
$H^{n-k}\left(G/H, \pi; K_{G/H}^{2-N}\right)$ coincides with the pairing 
between 
$H^{k}_{{\rm Lie}} (\l, H; C^\infty(G)_N)$ and $H^{n-k}_{{\rm Lie}}(\l, H; C^\infty(G)_{2-N})$
induced by $\nu$ (see $\S$\ref{relative}).

\subsection{Modular vector fields of $(G/H, \pi)$}\lb{modular-GH}
Assume again that $G/H$ is orientable and let $\mu$ be a fixed volume
form on $G/H$. 
Fix a non-zero $Y_0 \in \wedge^\top \h^0$, and let
$Y_{0}^{l}$ be the corresponding left invariant form on $G$. Write
$q^* \mu = \phi Y_{0}^{l}$, with 
 $\phi \in C^\infty(G)$ everywhere non-zero.
Let $\Lambda \in \wedge^2 \g$ be any element such that
$q(\Lambda) = \pi(eH) \in \wedge^2 T_{\ue} G/H \cong \wedge^2 \g/\h$, and
let $\pi_\Lambda = \Lambda^l+ \piG$ so that $q_* \pi_\Lambda = \pi$.
Recall that $\chi_\l \in \l^*, \chi_\g \in \g^*$ and $\chi_{\g^*} \in
\g$ are the adjoint characters of $\l, \g$ and $\g^*$ respectively. 
Write $x_0 = \chi_{\g^*} \in \g, \,\xi_0 = \chi_\g \in \g^*$, and let 
$x_\l$ be any element in $\g^*$ such that $x_\l(\xi) = \chi_\l(\Lambda \xi +
\xi)$ for $\xi \in \h^0$. Recall that for $x \in \g$ and $\xi \in \g$,
$x^l$ (resp. $x^r$ and $\xi^l$) is the left (resp. right) invariant vector field and
one form on $G$ with values $x$ and $\xi$ at $e \in G$.

\ble{modular-fields-GH}
Let the notation be as above. Let $X$ be the vector field on 
$G$ given by
\[
X = -\tpi_\Lambda (d \log|\phi|) + \frac{1}{2} \left(
x_{\l}^{l} + x_{0}^{r} + \tpi_\Lambda (\xi_{0}^{l}) \right).
\]
Then $q_*X$ is a well-defined vector field on $G/H$, and it is the
modular vector field of $\pi$ with respect to $\mu$.
\ele

\begin{proof} Let $\xi_1, \ldots, \xi_n$ be a basis of $\h^0$ such that
$\xi_1 \wedge \cdots \wedge \xi_n = Y_0$. 
Let $\alpha \in \Omega^1(G/H)$, and
write $q^* \alpha = \sum_{j=1}^{n} f_{\alpha, j} \xi_{j}^{l} \in \Omega^1(G)$. 
As in the proof of \thref{GH-canonical-rep},
\begin{eqnarray*}
q^* D_\alpha \mu & = & -\sum_{j=1}^{n} f_{\alpha, j} ((\xi_{j}^{l}, \, 
\tpi_\Lambda(d\log|\phi|) + (\piG, \, d \xi_{j}^{l}))Y_0^\l \\
& & + \frac{1}{2}\sum_{j=1}^{n} f_{\alpha, j} (\chi_\l(\Lambda \xi_j + \xi_j) -\chi_\g(\Lambda \xi_j)
+ \chi_{\g^*}(\xi_j))Y_0^l\\
& = & \left(q^* \alpha, \,\, -\tpi_\Lambda (d \log|\phi|)- F_0 
+ \frac{1}{2}(x_{\l}^{l} + (\Lambda \xi_0)^l + x_{0}^{l})\right),
\end{eqnarray*}
where $F_0$ is the vector field on $G$ such that $(F_0, \xi^l) = 
(\piG, d\xi^l)$ for all $\xi \in \g^*$. It is shown in 
Proposition 4.7 of \cite{e-lu-we} that 
$F_0 = \frac{1}{2} (x_{0}^{l} - x_{0}^{r} -\tilde{\pi}_{\scriptscriptstyle G}(\xi_{0}^{l})).$
Thus we have
\[
q^* D_\alpha \mu = (q^* \alpha, \; X).
\]
It follows that
$q_*X$ is a well-defined vector field on $G/H$ and it is the modular
vector field of $\pi$ with respect to $\mu$.
\end{proof}
 
\bre{hh0}
Note that if $\mu$ is a $G$-invariant volume form on $G/H$, the modular vector field
of $\pi$ with respect to $\mu$ is 
\[
q_* X = \frac{1}{2} \left(
x_{\l}^{l} + x_{0}^{r} + \tpi_\Lambda (\xi_{0}^{l}) \right).
\]
This formula for the special case when $\h^0$ is an ideal of $\g^*$ 
has been obtained in \cite{e-lu-we}.
\ere

\sectionnew{A Poisson groupoid over $(G/H, \pi)$}\lb{sym-oid}
When $H$ is a Poisson 
Lie subgroup of $(G, \piG)$ and $\pi = q_* \piG$, where $q: G \to G/H$ is the
projection, a symplectic groupoid of $(G/H, \pi)$ was constructed in
\cite{xuping-oid} (under the additional assumption that 
$(G, \piG)$ is complete). In this section, let $(G/H, \pi)$ be
an arbitrary Poisson homogeneous space of $(G, \piG)$ with the
Drinfeld Lagrangian subalgebra $\l = \l_{\pi(eH)}$.
We assume that
$G$ is a closed subgroup of a connected Lie group $D$ 
with Lie algebra  $\d$, $H = G \cap L$, where $L$ is the connected subgroup 
of $D$ with Lie algebra $\l$, and that  the infinitesimal
action $\lambda$ of $\l$ on $G$ in \eqref{lambda} integrates to
an action of   $L$ on $G$. We will
show that the associated space $\Gamma=G \times_H (L/H)$ is a Poisson groupoid over $(G/H, \pi)$.
We also give conditions for $\Gamma$ to be symplectic.
The Poisson structure on $\Gamma$ is obtained from reduction of
a quasi-Poisson manifold by an action of a quasi-Poisson Lie group
\cite{anton-yvette}.

\subsection{The quasi-Poisson Lie group $(G, \piGL, \varphi)$}\lb{subsec-quasi-G}
Let $(G, \piG)$ be a Poisson Lie group  corresponding to  
Manin triple $(\d, \g, \g^*)$. Then any $\Lambda \in \wedge^2 \g$ (not necessarily 
related to any Poisson homogeneous space of $(G, \piG)$ as in $\S$\ref{subsec-drinfi})
can be used to twist   the Manin triple $(\d, \g, \g^*)$
 to a Manin quasi-triple $(\d, \g, \g^\prime)$
\cite{anton-yvette}, where
\begin{equation}\lb{eq-g-prime}
\g^\prime = \{\Lambda \xi + \xi\,|\, \xi \in \g^*\},
\end{equation}
and thus defines a quasi-Poisson Lie group structure on $G$.
More precisely, let $p_1: \d \to \g: x + \xi \mapsto x - \Lambda \xi$, where
$x \in \g$ and $\xi \in \g^*$, be the projection from $\d = \g + \g^\prime$ 
to $\g$ along $\g^\prime$, and
define
$\varphi \in \wedge^3 \g$ by
\[
\varphi(\xi\wedge \eta \wedge \zeta) = \la p_1[\Lambda \xi + \xi, \, \Lambda \eta + \eta], \, 
\Lambda \zeta + \zeta \ra, \hs \xi, \eta, \zeta \in \g^*.
\]
It is straightforward to check that, for any $\xi, \eta, \zeta \in \g^*$,
\begin{align*}
\varphi(\xi\wedge \eta \wedge \zeta)&= \la [\Lambda \xi, \Lambda \eta], \, \zeta \ra + 
\la [\Lambda \eta, \Lambda \zeta], \, \xi \ra + \la [\Lambda \zeta, \Lambda \xi], \, \eta \ra \\
& \; \;\;+ \la \Lambda \xi, \, [\eta, \zeta] \ra + \la \Lambda \eta, \, [\zeta, \xi] \ra + 
\la \Lambda \zeta, \, [\xi, \eta] \ra.
\end{align*}
In fact $\varphi = \frac{1}{2}[\Lambda, \Lambda] + \delta (\Lambda).$ 
Let $\Lambda^l$ and $\Lambda^r$ be respectively the left and right invariant bi-vector fields on 
$G$ with
value $\Lambda$ at $e$, and define 
\begin{equation}\lb{eq-piGL}
\piGL = \Lambda^l - \Lambda^r + \piG,
\end{equation}

\ble{le-piGL} \cite{anton-yvette}
$(G, \piGL, \varphi)$ is a quasi-Poisson Lie group corresponding to the
Manin quasi-triple $(\d, \g, \g^\prime)$ in the sense that
$\piGL$ is multiplicative,
\[
\frac{1}{2}[\piGL, \, \piGL] = \varphi^l - \varphi^r, \hs
\mbox{and} \hs [\piGL, \, \varphi^l]=0,
\]
where $\varphi^r$ (resp. $\varphi^l$) is the right (resp. left)
invariant tri-vector field on $G$ with value $\varphi$ at $e$.
\ele

Recall from \cite{anton-yvette} that a  (right)
quasi-Poisson action of $(G, \piGL, \varphi)$
on a manifold $P$ with a bi-vector field $\piP$ is a right action 
$\rho: P \times G \to P$ of $G$ on $P$ 
such that

1) $[\piP, \piP] = 2\rho_\varphi$ and

2) $\rho: (P, \piP) \times (G, \piGL) \to (P, \piP)$ is a bi-vector map,

\noindent
where $\rho: x \mapsto \rho_x$ also denotes the Lie algebra homomorphism 
$ \g \to \V^1(P)$
given by
\begin{equation}\lb{eq-rho-x}
\rho_x(p) = \frac{d}{dt}|_{t=0} p \exp (tx), \hs x \in \g, \, p \in P
\end{equation}
as well as its multi-linear extension $\wedge^k \g \to \V^k(P): X \mapsto
\rho_X$ for any 
integer $k \geq 1$. Left quasi-Poisson actions of $(G, \piGL, \varphi)$ 
are similarly defined. 

\bex{ex-G-DG}
Let $\pi_\Lambda = \Lambda^l + \piG$. It is easy to see that the   action of
$(G, \piGL, \varphi)$ on $(G, \pi_\Lambda)$ by right multiplication is
a right quasi-Poisson action.  For another example,
assume that $D$ is a connected Lie group with Lie algebra $\d$ and  that
$G$ is a closed subgroup of $D$. For $d \in D$, let $\ud = dG \in D/G$, and for
$x + \xi \in \d$ with $x \in \g$ and $\xi \in \g^*$, let $\sigma_{x + \xi}$ be the
vector field on $D/G$ given by
\begin{equation}\lb{eq-sigma-x-xi}
\sigma_{x + \xi}(\ud) = \frac{d}{dt}|_{t=0} \exp (t(x + \xi)) \ud \in 
T_{\ud}(D/G), \hs d \in D.
\end{equation}
Let $\sigma: \wedge^k \d \to \V^k(D/G): X \mapsto \sigma_X$ also denote
the multi-linear extension of $\sigma$. Let $\{x_i\}_{n=1}^{n}$ be a basis of $\g$ 
and let $\{\xi_i\}_{i=1}^{n}$ be its dual basis of $\g^*$.
Define the bi-vector fields $\piDG$ and $\piDGL$ on  $D/G$  respectively
by 
\begin{equation}\lb{eq-piDGL}
\piDG =\frac{1}{2} \sum_{i} \sigma_{\xi_i} \wedge \sigma_{x_i} \hs
\mbox{and} \hs
\piDGL = \frac{1}{2} \sum_{i} \sigma_{\Lambda \xi_i + \xi_i} \wedge \sigma_{x_i}
=\piDG -\sigma_\Lambda.
\end{equation}
Then  \cite{anton-yvette}  $\piDG$ is Poisson and  the action
\[
 (G, \piGL, \varphi) \times (D/G, \, \piDGL) \lra(D/G, \, \piDGL):\;   
(g, \ud) \longmapsto g \ud, \hspace{.1in} g \in G, d \in D,
\]
is a left quasi-Poisson action of $(G, \piGL, \varphi)$. In 
particular,
\begin{equation}\lb{eq-piDGL-piDGL}
[\piDGL, \; \piDGL] = -2 \sigma_\varphi.
\end{equation}
Moreover, let $\delta_{\g^\prime}: \g^\prime \to \wedge^2 \g^\prime$ be
defined by
\[
\la \delta_{\g^\prime} (\Lambda \xi + \xi), \; x \wedge y \ra = 
\la \Lambda \xi + \xi, \, [x, y] \ra = \la \xi, \, [x, y]\ra, 
\hs \xi \in \g^*, \, x, y \in \g.
\]
Then one can check that
\begin{equation}\lb{xi-piDGL}
[\sigma_{\Lambda \xi + \xi}, \; \piDGL ] = -\sigma_{\delta_{\g^\prime}(\Lambda \xi + \xi)} + 
\sigma_{\iota_\xi \varphi}, \hs \forall \xi \in \g^*.
\end{equation}
\eex

\subsection{The bivector field $\piP$ on $P = G \times (D/G)$}\lb{subsec-P}
Let the assumptions be as in $\S$\ref{subsec-quasi-G}. In particular, assume that
$D$ is a connected Lie group with Lie algebra $\d$ and that $G$ is a closed subgroup of $D$.
Let $P = G \times (D/G)$.
For any integer $k \geq 1$ and for
a $k$-vector field $V$ on $G$, let  $(V, 0)$ be the corresponding $k$-vector field on
$P$. Similarly a $k$-vector field $U$ on $D/G$ gives rise to the $k$-vector field
$(0, U)$ on $P$. For $x \in \g$, recall that $x^l$ is the left invariant vector field on
$G$ with $x^l(e) = x$.
Define the bi-vector field $\piP$ on $P$ by
\begin{equation}\lb{eq-piP}
\piP = (\pi_\Lambda, 0) - (0, \piDGL) + \sum_{i=1}^{n} (0, \sigma_{\Lambda \xi_i + \xi_i})
\wedge (x_i^l, 0).
\end{equation}

\ble{le-piP}
The right action 
\[
\rho: \; (P, \piP) \times (G, \piGL, \varphi) \lra (P, \piP): \;  
(g, \ud) \cdot g_1 = (gg_1, \; g_1^{-1}\ud), \hs g, g_1 \in G, d \in D,
\]
is a quasi-Poisson action of $(G, \piGL, \varphi)$.
\ele

\begin{proof} To show that 
$[\piP, \piP] = 2\rho_\varphi$,
let $\varphi = \sum_k a_k \wedge b_k \wedge c_k$, where $a_k, b_k, c_k \in \g$,
and let
\begin{align*}
\rho_\varphi^\prime =& \sum_k \left(
(0, \sigma_{a_k}) \wedge (b_k^l \wedge c_k^l, 0) + 
(0, \sigma_{b_k}) \wedge (c_k^l \wedge a_k^l, 0) +
(0, \sigma_{c_k}) \wedge (a_k^l \wedge b_k^l, 0)
\right)\\
\rho_\varphi^{\prime\prime} = & \sum_k \left(
(a_k^l, 0) \wedge (0, \sigma_{b_k \wedge c_k}) + 
(b_k^l, 0) \wedge (0, \sigma_{c_k \wedge a_k}) + 
(c_k^l, 0) \wedge (0, \sigma_{a_k \wedge b_k})\right).
\end{align*}
It is easy to see that $\rho_\varphi = (\varphi^l, 0) - (0, \sigma_{\varphi})
-\rho_\varphi^\prime + \rho_\varphi^{\prime\prime}$. On the other hand, 
let $\pi_0 = \sum_{i=1}^{n} (0, \sigma_{\Lambda \xi_i + \xi_i})
\wedge (x_i^l, 0),$  so that $\piP = (\pi_\Lambda, 0) - (0, \piDGL) + \pi_0$, and
\begin{align*}
[\piP, \; \piP] &= ([\pi_\Lambda, \pi_\Lambda], 0) + (0, [\piDGL, \piDGL])
+2[\pi_0, \; (\pi_\Lambda, 0) -(0, \piDGL) ] + [\pi_0, \; \pi_0]\\
&=2(\varphi^l, 0) - 2(0, \sigma_\varphi) +2[\pi_0, \; (\pi_\Lambda, 0) -(0, \piDGL) ] + [\pi_0, \; \pi_0].
\end{align*}
It is easy to see that
$[\pi_0, \pi_0] = \pi_1 + \pi_2$, where
\begin{align*}
\pi_1 & = -\sum_{i, j = 1}^{n} (0, \sigma_{[\Lambda \xi_i +\xi_i, \Lambda \xi_j +\xi_j]})
\wedge ( x_i^l \wedge x_j^l, 0),\\
\pi_2 &= \sum_{i, j = 1}^n ([x_i, x_j]^l, 0) \wedge 
(0,  \sigma_{(\Lambda \xi_i +\xi_i) \wedge \Lambda \xi_j +\xi_j)}).
\end{align*} 
Thus  $[\piP, \piP] = 2(\varphi^l, 0) - 2(0, \sigma_\varphi)
+2[\pi_0, (\pi_\Lambda, 0)] + \pi_1 -2[\pi_0, (0, \piDGL)]  + \pi_2.$
Now
\[
2[\pi_0, (\pi_\Lambda, 0)] = 2\sum_{i=1}^n (0, \sigma_{\Lambda \xi_i + \xi_i}) 
\wedge (([x_i, \Lambda] +\delta(x_i))^l, 0).
\]
Recall that $p_1: \d \to \g$ is the projection along $\g^\prime$. Let $p^\prime: \d
\to \g^\prime$ be the projection along $\g$. It is easy to check that
\begin{align*}
\sum_{i,j=1}^{n} p^\prime[\Lambda \xi_i + \xi_i, \;\Lambda \xi_j + \xi_j] \otimes x_i \wedge x_j & =
2 \sum_{i=1}^n (\Lambda \xi_i + \xi_i) \otimes ([x_i, \Lambda] +\delta(x_i))\\
\sum_{i,j=1}^{n} p_1 [\Lambda \xi_i + \xi_i,\; \Lambda \xi_j + \xi_j] \otimes x_i \wedge x_j & =
2\tilde{\varphi},
\end{align*}
where $\tilde{\varphi} = \sum_k (a_k \otimes b_k \wedge c_k + b_k \otimes c_k \wedge a_k + 
c_k \otimes a_k \wedge b_k)$. Thus $2[\pi_0, (\pi_\Lambda, 0)] + \pi_1 = -2\rho_\varphi^\prime$.
Similarly, by 
\eqref{xi-piDGL}, 
\[
[\pi_0, \piDGL] = \sum_{i=1}^{n} (x_i^l, 0) \wedge (0, -[\sigma_{\Lambda \xi_i + \xi_i}, \, 
\piDGL])
= \sum_{i=1}^{n} (x_i^l, 0) \wedge (0, \sigma_{\delta_{\g^\prime}(\Lambda \xi_i + \xi_i)} - 
\sigma_{\iota_{\xi_i} \varphi}).
\]
It is easy to check that $\sum_{i=1}^{n} x_i \otimes \iota_{\xi_i} \varphi = \tilde{\varphi}$
and that
\[
2 \sum_{i=1}^{n} x_i \otimes \delta_{\g^\prime}(\Lambda \xi_i + \xi_i) = 
\sum_{i,j=1}^{n} [x_i, x_j] \otimes (\Lambda \xi_i + \xi_i) \wedge (\Lambda \xi_j + \xi_j).
\]
Thus $ -2[\pi_0, \piDGL]  +\pi_2 =2\rho_\varphi^{\prime\prime}.$ Hence
$[\piP, \piP] = 2\rho_\varphi$.
The proof that $\rho$ is a bi-vector map is straightforward and we omit the 
details.
\end{proof}

We now study when $\piP$ on $P = G \times (D/G)$ is nondegenerate. For
$d \in D$, the linear map
$\d \to T_{\ud}(D/G):  x + \xi \mapsto \sigma_{x + \xi}(\ud)$, where $x \in \g$ and $\xi \in \g^*$,
induces an isomorphism $\d/\Ad_d \g \to T_{\ud}(D/G)$. For $y \in \g$ and $\eta \in \g^*$, 
let $\alpha_{y + \eta}(\ud) \in T^*_{\ud}(D/G)$ be such that 
\begin{equation}\lb{eq-alpha-y-eta}
(\alpha_{y + \eta}(\ud), \, \sigma_{x + \xi}(\ud) ) = \la y + \eta, \; x + \xi \ra, \hs x \in \g, \xi \in 
\g^*.
\end{equation}
Then we have the isomorphism
\begin{equation}\lb{eq-alpha-iso}
\Ad_d \g \lra T^*_{\ud}(D/G): \; \; y + \eta \longmapsto \alpha_{y + \eta}(\ud), \hs y \in \g, \, 
\eta \in \g^*, \, y + \eta \in \Ad_d \g.
\end{equation}
Note that when $y + \eta \in \Ad_d \g$, $\sigma_{y + \eta}(\ud) = 0$, so 
$\sigma_y(\ud) = -\sigma_\eta(\ud)$. The proof of the first identity in
the following \leref{le-tilde-piP} is
straightforward and is omitted. The second identity follows from \eqref{lambda-1}.

\ble{le-tilde-piP}
For $g \in G, d \in D, \xi \in \g^*$ and $y + \eta \in \Ad_d \g$ with $y \in \g$ and $\eta 
\in \g^*$,
\begin{align*}
\tilde{\pi}_{\scriptscriptstyle P}(l_{g^{-1}}^{*}\xi, \, \alpha_{y + \eta}(\ud) ) & =
\left( \tilde{\pi}_{\scriptscriptstyle G} (l_{g^{-1}}^{*}\xi)+
(l_g)_* (y + \Lambda \xi -\Lambda \eta),\;\;
\sigma_{\Lambda \eta + \eta - \Lambda \xi - \xi}(\ud)\right)\\
& = \left(\lambda_{y - \Lambda \eta + \Lambda \xi + \xi}(g), \;\; 
-\sigma_{y - \Lambda \eta + \Lambda \xi + \xi}(\ud)\right)
\end{align*}
\ele

\ble{le-piP-nondegenerate}
The bi-vector field $\piP$ on $P = G \times (D/G)$ is nondegenerate at $(g, \ud)$, where
$g \in G$ and $d \in D$, if
\begin{equation}\lb{eq-gg}
\g^\prime \cap \Ad_d \g = 0 \hs \mbox{and} \hs \g^* \cap \Ad_{gd}\g = 0.
\end{equation}
In particular, $\piP$ is nondegenerate 
at $(g, \underline{e})$ for any $g \in G$, where $e \in D$ is the identity.
\ele

\begin{proof} Assume that \eqref{eq-gg} holds at $(g, \ud) \in P$. Suppose that
$\xi \in \g^*$ and $y + \eta \in \Ad_d \g$ are such that 
$\tilde{\pi}_{\scriptscriptstyle P}(l_{g^{-1}}^{*}\xi, \, \alpha_{y + \eta}(\ud) )=0$.
Then $\sigma_{\Lambda \eta + \eta - \Lambda \xi - \xi}(\ud)=0$
by \leref{le-tilde-piP}, so
$\Lambda \eta + \eta - \Lambda \xi - \xi \in \g^\prime \cap \Ad_d \g=0$. Thus $\xi = \eta$. 
By \eqref{Gond}, $\tilde{\pi}_{\scriptscriptstyle G} (l_{g^{-1}}^{*}\xi)
=(r_g)_* p_\g \Ad_g \xi$, where $p_\g: \d \to \g$ is the projection along $\g^*$.
Thus \leref{le-tilde-piP} implies that $p_\g \Ad_g (y + \xi) = 0$, so 
$\Ad_g (y + \xi) \in \g^* \cap \Ad_{gd}\g = 0$. Thus $y = 0$ and $\xi = \eta = 0$.
\end{proof}

\bre{re-nondegenerate-piP}
Let $N(\g^*)$ be the normalizer subgroup of $\g^*$ in $D$. Suppose that
$D = N(\g^*)G$ and that $\Lambda = 0$ (so $\pi(eH) = 0)$. Then 
\eqref{eq-gg} holds for all $(g, d) \in G \times D$, and $\piP$ is 
nondegenerate everywhere on $P$. See \exref{ex-semi-simple} for an example.
\ere

\subsection{The Poisson structure $\Pi$ on $G \times_H (L/H)$}\lb{subsec-Pi-GHLH}
Let the notation be as in $\S$\ref{subsec-quasi-G} and  $\S$\ref{subsec-P}, but assume
now that $(G/H, \pi)$ is a Poisson homogeneous space of $(G, \pi)$ and that
$\Lambda \in \wedge^2 \g$ is such that
$q(\Lambda) = \pi(eH) \in \wedge^2 T_{eH}(G/H) \cong \wedge^2 (\g/\h)$, where
$q$ denotes both projections $G\to G/H$ and $\g \to \g/\h$. Let $(G, \piGL, \varphi)$ be the
quasi-Poisson Lie group defined using $\Lambda$ as in $\S$\ref{subsec-quasi-G}.

\ble{le-PH} Let $P$ be any manifold with a bi-vector field $\piP$.
Suppose that $\rho: (P, \piP) \times (G, \piGL, \varphi) \to (P, \piP)$ is a right
quasi-Poisson action of $(G, \piGL, \varphi)$ and  that $\rho$ restricts to a 
free and proper action of $H$. Let $j: P \to P/H$ be the projection. Then 
$j_* \piP$ is a well-defined Poisson structure on $P/H$.
\ele

\begin{proof} By 1) in \leref{le-Lambda}, $q_* \piGL(h) = 0$ for all $h \in H$. It
follows from the fact that $\rho$ is a bi-vector map that
$j_* \piP$ is  well-defined.
Since $\varphi \in \h \wedge \g \wedge \g$ by 2) of \leref{le-Lambda}, 
$[j_* \piP, \, j_* \piP] = j_*[\piP, \, \piP] = 2j_* \rho_\varphi = 0$, so
$j_* \piP$ is Poisson.
\end{proof}

We now state a lemma from linear algebra.

\ble{le-linear-algebra} Let $(V, \pi)$ be a Poisson vector space. Suppose that 
$U$ and $W$ are subspaces of $V$ such that $\tpi(U^0) \subset W \subset U$, where
$U^0 = \{\xi \in V^* \, |\, \xi|_U = 0\}$. Let 
$\phi: V \to V/W$ be the projection. Then $U/W$ is a Poisson subspace of
$(V/W, \phi(\pi))$.
\ele

The following \leref{le-QH} follows immediately from \leref{le-linear-algebra}.

\ble{le-QH} Let the notation be as in \leref{le-PH}. Suppose that $Q$ is an $H$-invariant 
submanifold of $P$ such that $\tilde{\pi}_{\scriptscriptstyle P}(T^0_qQ) \subset T_q(qH)$
for every $q \in Q$, where $T^0_qQ = \{\alpha \in T^*_qP | \, \alpha|_{T_qQ} = 0\}$ and
$qH$ is the $H$-orbit through $q$. Then $Q/H$ is a Poisson submanifold of
$(P/H, \, j_*\piP)$.
\ele

We now apply \leref{le-PH} to  $P=
G \times (D/G)$ as in $\S$\ref{subsec-P}, $\piP$ as in
\eqref{eq-piP}, and the action $\rho$ as in  \leref{le-piP}. 
Denote by $G \times_H (D/G)$
the quotient of $P$ by $H$ with the projection
$j: P \to G \times_H (D/G)$. By \leref{le-PH}, $j_*\piP$ 
is a well-defined Poisson structure on $G \times_H (D/G)$.
Set $[g, \ud] = j(g, \ud)$ for $g \in G$ and $d \in D$.

\bnota{nota-Pi}
The Poisson structure $j_* \piP$ on $G \times_H (D/G)$ will be denoted by $\Pi$.
\enota

Recall that $\l = \l_{\pi(eH)}$ is the Drinfeld Lie subalgebra of $\d$ 
associated to $\pi(eH)$. 
Let $L$ be the connected Lie subgroup of $D$ with Lie algebra $\l$ and assume that
$H = G \cap L$. Let $\O$ be the $L$-orbit in $D/G$ through $\underline{e} 
\in D/G$, where $e$ is the identity element of $D$.
 Identify $L/H$ with $\O$ and regard $G \times_H (L/H)$ as a submanifold of $G \times_H (D/G)$.

\ble{le-GHLH}
$G \times_H (L/H)$ is a Poisson submanifold of $(G \times_H (D/G), \; \Pi)$, and 
$\Pi$ is nondegenerate at $[g, \ud]$ for all $g \in G$ and $d \in L$ such that
\eqref{eq-gg} holds.
\ele

\begin{proof} Let $Q = G \times \O \subset P$. Then $Q$ is $H$-invariant. To
see that $Q/H$ is a Poisson submanifold of $(P/H, \Pi)$, it suffices, by 
\leref{le-QH}, to show that $\tilde{\pi}_{\scriptscriptstyle P}(T^0_qQ) \subset T_q(qH)$
for every $q = (g, \ud) \in Q$, where $g \in G$ and $d \in L$. Using the
isomorphism in \eqref{eq-alpha-iso}, 
$T^0_qQ = \{(0, \alpha_{y + \eta}(\ud)) \, | \, y + \eta \in \l \cap \Ad_d \g\}$, and by
\leref{le-tilde-piP}, 
\[
\tilde{\pi}_{\scriptscriptstyle P}(T^0_qQ) = \{((l_g)_* (y - \Lambda \eta), \; 
\sigma_{\Lambda \eta - y}(\ud))\, | \, y + \eta \in \l \cap \Ad_d \g\}.
\]
By \eqref{eq-l-Lambda}, $y + \eta \in \l$ implies that $y - \Lambda \eta\in \h$. Thus
$\tilde{\pi}_{\scriptscriptstyle P}(T^0_qQ) \subset T_q(qH)$.

By \leref{le-piP-nondegenerate}, $\piP$ is nondegenerate at $(g, \ud)$ for 
all $g \in G$ and $d \in D$ such that \eqref{eq-gg} holds. At such a point
$(g, \ud)$ where $d \in L$,  $\l \cap \Ad_d \g = \Ad_d (\l \cap \g) = \Ad_d \h$, and
the map $\l \cap \Ad_d \g \to \h: y + \eta \mapsto y - \Lambda \eta$ is an 
isomorphism, so $\tilde{\pi}_{\scriptscriptstyle P}(T^0_qQ)= T_q(qH)$.
It follows from a linear algebra argument that  
$\Pi$ is nondegenerate at $[g, \ud]$.
\end{proof}

\subsection{The Poisson groupoid $(G \times_H (L/H), \Pi)$}\lb{subsec-poi-gpoid}

Let the notation be as in $\S$\ref{subsec-Pi-GHLH}. Recall that
$\lambda: \d \to \V^1(G)$ is the infinitesimal action of $\d$ on $G$ 
given in \eqref{lambda}. Assume in addition that
the restriction of $\lambda$ to $\l$ integrates to a right action of $L$ on $G$, denoted by
$(g, l) \mapsto g^l$ for $g \in G$ and $l \in L$,  such that
$g^h = gh$ for $g \in G$ and $h \in H$. Then $G$ is a $(\d, L)$-space (see 
\deref{lHspace}).

Let $\Gamma = G \times_H (L/H)$.
It is straightforward to show (we omit the proof) that 
the following is a groupoid structure on $\Gamma$ over $G/H$: for $g, g_1, g_2 \in G$
and $l, l_1, l_2 \in L$,

1) source map $s: \Gamma \to G/H: \, [g,\, lH] \mapsto gH$;

2) target map $t: \Gamma \to G/H: \, [g, \,lH] \mapsto g^lH$;

3) multiplication $\cdot_\Gamma$: $[g_1,\, l_1H]\cdot_{\scriptscriptstyle \Gamma}  [g_2,\, l_2H] =
[g_1, \;l_1 h l_2H]$ when $g_1^{l_1}H = g_2H$, where $h = (g_1^{l_1})^{-1}g_2$;

4) inverse $\tau: \Gamma \to \Gamma: \, [g, \,lH] \mapsto [g^l, \,l^{-1}H]$; 

5) identity section $\epsilon: G/H \to \Gamma: \, gH \mapsto [g, \,eH]$.

\bth{th-poi-grpoid}
With the groupoid structure described above,  $(G \times_H (L/H), \Pi)$
is a Poisson groupoid over $(G/H, \pi)$.
\eth

The proof of \thref{th-poi-grpoid} will be given in $\S$\ref{subsec-proof}. 

\bre{re-symplectic}
Assume that $\pi(eH) = 0$, so that we can take $\Lambda = 0$. Recall that
$G^*$ is the 
connected subgroup of $D$ with Lie algebra $\g^*$. Assume further that the map
$G^* \times G \to D: (u, g) \mapsto ug$ is a diffeomorphism. Identify $G$ 
with $G^* \backslash D$. Then 
the restriction to $L$ of the right action of $D$ on $G \cong G^* \backslash D$ integrates
the infinitesimal action $\lambda$ of $\l$ on $G$. By \reref{re-nondegenerate-piP} and
\leref{le-GHLH},  
$\Pi$ is nondegenerate everywhere on
$G \times_H (L/H)$. Thus $(G \times_H (L/H), \Pi)$
is a symplectic groupoid over $(G/H, \pi)$. Note that in this case, the bi-vector field
$\piP$ on $P \cong D$ is Poisson by \leref{le-piP} and everywhere nondegenerate by 
\leref{le-piP-nondegenerate}. Moreover, by \leref{le-h0}, $\h^0$ is a subalgebra of $\g^*$. 
Let $H^0$ be the connected subgroup of $G^*$ with Lie algebra $\h^0$. Then $L = HH^0$ is
a coisotropic submanifold $(D, \piP)$ and $L/H \cong H^0$. Our construction of the
symplectic structure $\Pi$ on 
$G \times_H (L/H) \cong G \times_H H^0$
is a special case of coisotropic reduction for symplectic manifolds. 
In the special case when $H$ is a Poisson
subgroup of $(G, \piG)$ and when $\pi = q_* \piG$, this construction was carried out 
in \cite{xuping} 
\ere

\bex{ex-semi-simple}
Let $G$ be a connected and simply connected Lie group and let $X$ be the variety of
Borel subgroups of $G$. Let $G_0$ be a real form of $G$ and  $K$  a compact real form
of $G$ such that $K_0:=G \cap K$ is a maximal compact subgroup of $G_0$. Choose
an Iwasawa decomposition $G = KAN$ of $G$ such that the Borel subgroup $B$ of $G$  
containing $AN$ 
lies in the unique closed $G_0$-orbit in $X$. This 
choice of $B$ gives rise to a Poisson Lie group
$(K, \pi_{{\scriptscriptstyle K}})$ with $AN$ as a dual Poisson Lie group. Although
$K_0$ is not a Poisson Lie subgroup of
$(K, \pi_{{\scriptscriptstyle K}})$, it is shown in \cite{foth-lu:symmetric} that   
the projection $\pi$ of $\pi_{{\scriptscriptstyle K}}$ is a 
well-defined 
Poisson structure  on $K/K_0$, making
$(K/K_0, \pi)$ a Poisson homogeneous space of $(K, \pi_{{\scriptscriptstyle K}})$,
and the Drinfeld Lagrangian subalgebra associated to $\pi(eK_0)$ is $\g_0$, the Lie algebra
of $G_0$.
Let $T = K \cap B$, a
maximal torus of $K$.  
The set of $T$-orbits of symplectic leaves of 
$\pi$ in $K/K_0$ is shown in \cite{foth-lu:symmetric} to be
in one to one correspondence with the set of $G_0$-orbits in $X$.
Due to the importance in representation theory of $G_0$-orbits in $X$, 
the Poisson geometrical properties of 
$(K/K_0, \pi)$ are worth further study.
Since $\pi(eK_0) = 0$ and since $G = KAN = ANK$,  
the conditions in \reref{re-nondegenerate-piP}
are satisfied. By \reref{re-symplectic} and \thref{th-poi-grpoid},
$K \times_{K_0} (G_0/K_0)$ has the structure of a symplectic
groupoid over $K/K_0$. More details of this example, in particular, the generalized 
Poisson cohomology of $(K/K_0, \pi)$, will be studied in a future paper.
\eex

\subsection{Proof of \thref{th-poi-grpoid}}\lb{subsec-proof}

Let the assumptions be as in $\S$\ref{subsec-poi-gpoid}. We  need two lemmas.
Recall that
$G^*$ is the connected subgroup of $D$ with Lie algebra $\g^*$. 

\ble{le-in-Gs}
For any $g \in G$ and $l \in L$, $g^l l^{-1} g^{-1} \in G^*$.
\ele

\begin{proof}
Fix $g \in G$ and $l \in L$. To avoid confusion with the notation set
up in $\S$\ref{notation} for left and right translations  on $G$, 
if $v  \in T_gG$, we let $g^{-1} v  \in \g$ and
$v  g^{-1} \in \g$ be the left and right translation of $v $ by $g^{-1}$.

Let $l(t)$ be a smooth curve in $L$ such that $l(0) = e$ and $l(1) = l$, and let
$u(t) = g^{l(t)} l(t)^{-1} g^{-1} \in D$. Let $u^\prime(t) \in T_{u(t)}D$ and
$l^\prime(t) \in T_{l(t)}L$ be respectively the derivatives of $u(t)$ and $l(t)$ at $t$.
Let $x(t) = l(t)^{-1} l^\prime(t) \in l$.
Then, for every $t$,
\[
u^\prime(t) = \lambda_{x(t)}(g^{l(t)}) l(t)^{-1} g^{-1} - g^{l(t)} x(t) l(t)^{-1} g^{-1}
\]
so  by \eqref{lambda-1}, 
$u^\prime(t) u(t)^{-1} = -\Ad_{g^{l(t)}} x(t) +p_\g \Ad_{g^{l(t)}} x(t) = 
-p_{\g^*} \Ad_{g^{l(t)}} x(t) \in \g^*,$
where $p_\g$ and $\p_{\g^*}$ are projections from $\d$ to  $\g$ and $\g^*$ with 
respect to the decomposition $\d = \g + \g^*$. It follows from $u(0) = e$ that
$u(t) \in G^*$ for all $t$. In particular, $g^l l^{-1} g^{-1} 
=u(1) \in G^*$.
\end{proof}

The following \leref{le-xi-h} is equivalent to 1) in \leref{le-Lambda}, and we omit its proof.

\ble{le-xi-h}
One has
$\Ad_h \Lambda \xi + p_\g \Ad_h \xi - \Lambda \Ad_{h^{-1}}^* \xi \in \h$
for all $\xi \in \h^0$ and $h \in H$.
\ele
 
We can now start the 
proof of \thref{th-poi-grpoid}.

Let $\G_\Gamma = \{(\gamma_1, \gamma_2, \gamma_3) \in \Gamma \times \Gamma \times \Gamma \,
| \, t(\gamma_1) = s(\gamma_2), \, \gamma_3 = \gamma_1 
\cdot_{\scriptscriptstyle \Gamma} \gamma_2\}.$
By the definition of Poisson groupoids \cite{we:spoid}, we need to show that $\G_\Gamma$
is coisotropic in $\Gamma \times \Gamma \times \Gamma$ with the Poisson 
structure $\Pi \oplus \Pi \oplus (-\Pi)$. Let $(P, \piP)$ be as in
$\S$\ref{subsec-P} and recall that $j: P \to P/H$ is the
projection.
Since $(\Gamma, \Pi)$ is a Poisson submanifold of $(P/H, \Pi)$, and since
$j: (P, \piP) \to (P/H, \Pi)$ is a bi-vector map, it is enough 
\cite[Corollary 2.2.5]{we:coisotropic} to show that 
$\G_P:=(j \times j \times j)^{-1}(\G_\Gamma)$ is coisotropic in $(P \times P \times
P, \; \piP \oplus \piP \oplus (-\piP))$. 
Recall that $\O$ is the $L$-orbit in $D/G$ through $\ue\in D/G$ and that $Q = 
G \times \O \subset P$. Identify$L/H$ with $\O$ by identifying 
$lH \in L/H$ with  $\underline{l} = lG \in D/G$ for $l \in L$. Then
\begin{align*}
\G_P = &\{\left((g_1,\, l_1H),\, (g_2, \,l_2H), \,
(g_1h_3,\, h_3^{-1} l_1 hl_2H)\right) \, |\,\\
& \; \; \; g_1, g_2, g_2 \in G, \, l_1, l_2 \in L,  
\, h_3 \in H, \, g_1^{l_1}H = g_2H, \,
h = (g_1^{l_1})^{-1}g_2\} \subset Q \times Q \times Q.
\end{align*}
We will first describe the tangent bundle of $\G_P$ and then the 
co-normal bundle of $\G_P$ in $P \times P \times P$.

Let $\G_2 = \{\left((g_1,\, l_1H),\, (g_2, \,l_2H)\right) 
\, |\, g_1, g_2 \in G, \, l_1, l_2 \in L,\,
g_1^{l_1}H = g_2H\} \subset Q \times Q$. 
We now compute
$T_{(q_1, q_2)}\G_2$ for 
$(q_1, q_2) = ((g_1,\, l_1H), \, (g_2, \,l_2H)) \in \G_2$. 
Define $\tilde{t}, \tilde{s}: Q \to G/H$  by
$\tilde{t}(g, lH) = g^lH$ and $\tilde{s}(g, lH) = gH$ for $g \in G$ and $l \in L$.
Then 
\[
T_{(q_1, q_2)}\G_2 = \{(v_1, v_2)\, | \, v_1 \in T_{q_1}Q, \, v_2 \in T_{q_2}Q, \, 
\tilde{t}_*(v_1) = \tilde{s}_*(v_2)\}.
\]
Recall that $\sigma: \d \to \V^1(D/G)$ is given in \eqref{eq-sigma-x-xi}.
Let $\kappa: \g \to \V^1(G/H): x \to \kappa_x$ be the Lie 
algebra anti-homomorphism given by
\begin{equation}\lb{eq-kappa}
\kappa_x(gH) = \frac{d}{dt}|_{t = 0} \exp (tx) gH, \hs x \in \g, \,  g \in G.
\end{equation}
For $x, z \in \g$, $\zeta \in \g^*$ with $z + \zeta \in \l$, and $q = (g, lH) \in Q$, let
\[
v_{x, \,z + \zeta}(q) = \left( (l_g)_*x, \;\,  \sigma_{z + \zeta}(lH)\right) \in T_qQ.
\]
(Recall from $\S$\ref{notation} that the ``$l$" in $l_g$ denotes the left translation by $g$.
This is not to be confused with an element in $L$.) 
Recall that $p_\g: \d \to \g$ is the projection along $\g^*$.
Using the fact that $G$ is a $(\d, L)$-space via the infinitesimal action 
$\lambda$ of $\d$ and the action of $L$, one sees that 
\[
\tilde{t}_* v_{x_1, \,z_1 + \zeta_1}(q_1) = \kappa_{p_\g 
\Ad_{\scriptscriptstyle{g_1^{l_1} l_1^{-1}}}(x_1 + z_1 + \zeta_1)} (g_1^{l_1}H), \hs \mbox{for} \; 
x_1 \in \g, \,z_1 +\zeta_1 \in \l.
\]
Since $\tilde{s}_*v_{x_2, \,z_2 + \zeta_2}(q_2) = \kappa_{\Ad_{g_2} x_2}(g_2H)$ for 
$x_2\in \g, \,z_2 +\zeta_2 \in \l$, we get
\begin{align}\lb{eq-TG2}
T_{(q_1, q_2)}\G_2 & = \left\{\left(v_{x_1, \,z_1 + \zeta_1}(q_1), \;\,
v_{x_2, \,z_2 + \zeta_2}(q_2)\right) \, | \, x_1, x_2 \in \g, \; z_1 + \zeta_1, z_2 + \zeta_2 \in
\l,\right.\\
\nonumber
& \; \; \;  \;\;\;x_2 = \bar{x}_2+  \Ad_{g_2^{-1}} p_\g 
\Ad_{\scriptscriptstyle{g_1^{l_1} l_1^{-1}}}(x_1 + z_1 + \zeta_1) \, \,
\mbox{for some} \,\, \bar{x}_2 \in \h  \}.
\end{align}
Define 
\[
f: \; \G_2 \lra Q: \; ((g_1,\, l_1H), \, (g_2, \,l_2H)) \longmapsto 
(g_1, \, l_1 (g_1^{l_1})^{-1}g_2l_2H).
\]
Fix $g_i \in G, l_i \in L$, for $i = 1, 2$,
such that  
$(q_1, q_2) = ((g_1,\, l_1H), \, (g_2, \,l_2H)) \in \G_2$. Let 
$x_i \in \g$ and $z_i + \zeta_i \in \l$
be such that 
$(v_{x_1, \,z_1 + \zeta_1}(q_1),\;
v_{x_2, \,z_2 + \zeta_2}(q_2)) \in T_{(q_1, q_2)}Q$ as in \eqref{eq-TG2}. Let 
$g_1(t), l_1(t)$, and $l_2(t)$ be smooth curves in $G$ and $L$ respectively such that
$g_1(0) = g_1, \, g_1^\prime(0) = (l_{g_1})_* x_1$, and $ l_i(0) = l_i, \, 
l_i^\prime(0) = (r_{l_i})_* (z_i + \zeta_i)$ for $i = 1, 2$, where the superscript $\prime$ denotes
derivative at $0$. Let $g_2(t) = g_1(t)^{l_1(t)}
h \exp t\bar{x}_2$, where $h = (g_1^{l_1})^{-1} g_2 \in H$. It is easy to see that
$g^\prime_2(0) = (l_{g_2})_* x_2$. Let
\[
c(t) = \left( (g_1(t), \; l_1(t)H), \; (g_2(t), \;  l_2(t)H) \right) \in Q \times Q.
\]
Then $c(t) \in \G_2$ for all $t$, $c(0) = (q_1, q_2)$, and 
$c^\prime(0) = (v_{x_1, \,z_1 + \zeta_1}(q_1),\;\,
v_{x_2, \,z_2 + \zeta_2}(q_2))$. Since 
$f(c(t)) =(g_1(t), \; l_1(t)h \exp t\bar{x}_2 l_2(t)H)$, we have 
\[
f_* (v_{x_1, \,z_1 + \zeta_1}(q_1),\;\,
v_{x_2, \,z_2 + \zeta_2}(q_2)) = \frac{d}{dt}|_{t = 0}f(c(t)) 
= v_{x_1, \,z_3 + \zeta_3}(g_1, \, l_1hl_2H)),
\]
where  $z_3 \in \g$ and $\zeta_3 \in
\g^*$ are such that
\begin{equation}\lb{eq-zz3}
z_3 + \zeta_3 = 
z_1 + \zeta_1 + \Ad_{l_1h} (\bar{x}_2 + \Ad_{l_2} \bar{x}_3 +z_2 +\zeta_2) \hs \mbox{for 
some} \; \; \bar{x}_3 \in \h.
\end{equation}
Thus for $(q_1, q_2, q_3)\! = \!((g_1, l_1H),\, (g_2, l_2H), \,
(g_1h_3,\, h_3^{-1} l_1 hl_2H)) \in \G_P$, where $h_3 \in H$,
\begin{align*}
T_{(q_1, q_2, q_3)}\G_P &=  \left\{\left(v_{x_1, \,z_1 + \zeta_1}(q_1), \;\,
v_{x_2, \,z_2 + \zeta_2}(q_2), \; \, 
v_{x_3 + \Ad_{h_3}^{-1} x_1, \, -x_3 + \Ad_{h_3}^{-1} (z_3 + \zeta_3)}(q_3) \right)\, | \,
\right.\\
& \; \; \;  \;\;\; x_1, x_2 \in \g, \, x_3 \in \h, \; z_i \in \g, 
\zeta_i \in \g^*, \, z_i + \zeta_i
\in \l, \, i = 1, 2, 3,
\\
& \; \; \;  \;\;\;x_2 = \bar{x}_2+  \Ad_{g_2^{-1}} p_\g 
\Ad_{\scriptscriptstyle{g_1^{l_1} l_1^{-1}}}(x_1 + z_1 + \zeta_1) \, \,
\mbox{for some} \,\, \bar{x}_2 \in \h, \\
& \; \; \;  \;\; 
\left.z_3 + \zeta_3 = z_1 + \zeta_1 + \Ad_{l_1h} 
(\bar{x}_2 + \Ad_{l_2} \bar{x}_3+ z_2 +\zeta_2) \, \,
\mbox{for some} \,\, \bar{x}_3 \in \h\right\}.
\end{align*}
Let $T_{(q_1, q_2, q_3)}^0\G_P$
be the co-normal subspace of $T_{(q_1, q_2, q_3)}\G_P$ in 
$T_{(q_1, q_2, q_3)}^*(P \times P \times P)$.
Recall that for $\ud \in D/G$,  $\alpha_{y + \eta}(\ud) \in T^*_{\ud}(D/G)$ is given in
\eqref{eq-alpha-y-eta}. For  $y \in \g$, $\xi, \eta \in 
\g^*$,
and $q=(g, lH) \in Q$, let $\alpha_{\xi, \, y + \eta}(q) 
=(l_{g^{-1}}^* \xi, \, \alpha_{y + \eta}(lH))
\in T^*_qP$. Then  for 
$x, z \in \g$ and $ \zeta \in \g^*$ with 
$z + \zeta \in \l$,
\begin{equation}\lb{eq-alpha-v}
(\alpha_{\xi,\, y + \eta}(q), \; v_{x, \, z + \zeta}(q)) = 
(x, \xi) + \la y + \eta, \, z + \zeta\ra =
(x, \xi) + ( y, \zeta) + (z, \eta),
\end{equation}
and the map $\g^* \times \Ad_l \g \to T_q^*P: (\xi, \, y + \eta) \mapsto \alpha_{\xi, \, 
y + \eta}(q)$ is an isomorphism. Let $l_3 = h_3 l_1 h l_2 \in L$, where
$h = (g_1^{l_1})^{-1} g_2$. It follows 
from \eqref{eq-alpha-v} that $T_{(q_1, q_2, q_3)}^0\G_P$ consists of all
triples
\begin{equation}\lb{eq-alpha-triple}
(\alpha_{\xi_1, \, y_1 + \eta_1}(q_1), \;  \alpha_{\xi_2, \, y_2 + \eta_2}(q_2), \;
\alpha_{\xi_3, \, y_3 + \eta_3}(q_3)) \in T_{(q_1, q_2, q_3)}^*(P \times  P \times P), 
\end{equation}
where  $\xi_i, \eta_i \in \g^*, \,y_i \in \g$ and 
$y_i + \eta_i \in \Ad_{l_i} \g$ for $i = 1, 2, 3$, such that
\begin{align*}
0 &= (\xi_1, x_1) + (\xi_2, \, \Ad_{g_2^{-1}} p_\g \Ad_{g_1^{l_1}l_1^{-1}}x_1) + 
(\xi_3, \, \Ad_{h_3^{-1}} x_1) \\
&  +  \la y_1 + \eta_1, \, z_1 + \zeta_1 \ra + 
(\xi_2, \, \Ad_{g_2^{-1}} p_\g \Ad_{g_1^{l_1}l_1^{-1}}(z_1 + \zeta_1) ) 
+ \la y_3 + \eta_3, \, \Ad_{h_3^{-1}} (z_1 + \zeta_1) \ra \\
&   + (\xi_2, \bar{x}_2) + \la y_3 + \eta_3, \, \Ad_{h_3^{-1} l_1 h} \bar{x}_2 \ra \\
&   + \la  y_2 + \eta_2, \, z_2 + \zeta_2\ra 
+ \la y_3 + \eta_3, \, \Ad_{h_3^{-1} l_1 h} (z_2 + \zeta_2)\ra\\
&+ (\xi_3, x_3) -\la y_3 + \eta_3, \, x_3 \ra
\end{align*}
for all $x_1 \in \g, \, \bar{x}_2, x_3 \in \h, \, \z_1 + \zeta_1 \in \l$ and 
$z_2 + \zeta_2 \in \l$, which is equivalent to 
\begin{align}
\lb{eq-1}
&\xi_1 + \Ad_{l_1 (g_1^{l_1})^{-1}} \Ad_{g_2^{-1}}^* \xi_2 + \Ad_{h_3 }\xi_3 \in \g; \\
\lb{eq-2}
&y_1 + \eta_1 + \Ad_{l_1 (g_1^{l_1})^{-1}} \Ad_{g_2^{-1}}^* \xi_2 
+ \Ad_{h_3} (y_3 + \eta_3) \in \l;\\
\lb{eq-3}
&\xi_2 + \Ad_{h^{-1} l_1^{-1} h_3} (y_3 + \eta_3) \in \g +\l; \\
&\lb{eq-4} 
y_2 + \eta_2 +  \Ad_{h^{-1} l_1^{-1} h_3} (y_3 + \eta_3) \in \l; \\
\lb{eq-5} 
&\xi_3 - \eta_3 \in \h^0,
\end{align}
where recall that $\h^0 = \{\xi \in \g^* | (\xi, \h) = 0\}$.
Since $\g + \l = \g + \h^0$ by \eqref{eq-l-Lambda}, \eqref{eq-3} and \eqref{eq-4} imply that
$\eta_2 - \xi_2 \in \h^0$. Similarly, \eqref{eq-1}, \eqref{eq-2}, and
\eqref{eq-5} imply that $\eta_1 - \xi_1 \in \h^0$. Thus
\begin{equation}\lb{eq-in-h0}
\Lambda(\eta_i - \xi_i) + \eta_i - \xi_i \in \l \hs \mbox{for} \; \; i = 1, 2, 3.
\end{equation}

It remains to show that for any triple  in \eqref{eq-alpha-triple} satisfying  
\eqref{eq-1} - \eqref{eq-5},  
\begin{equation}\lb{eq-alpha-triple-ok}
\left( \tilde{\pi}_{\scriptscriptstyle P}(\alpha_{\xi_1, \, y_1 + \eta_1}(q_1)), \;
\tilde{\pi}_{\scriptscriptstyle P}(\alpha_{\xi_2, \, y_2 + \eta_1}(q_2)), \;
-\tilde{\pi}_{\scriptscriptstyle P}(\alpha_{\xi_3, \, y_3 + \eta_3}(q_3))\right)
\in T_{(q_1, q_2, q_3)} \G_P.
\end{equation}
Let $g_3 = g_1 h_3$. By \leref{le-tilde-piP} and \eqref{eq-in-h0}, 
\[
\tilde{\pi}_{\scriptscriptstyle P}(\alpha_{\xi_i, \, y_i + \eta_i}(q_i))
=v_{x_i,\, z_i + \zeta_i}(q_i) \in T_{q_i} Q
\]
for $i = 1, 2, 3$,  and 
\begin{equation}\lb{eq-i}
x_i = y_i + \Lambda \xi_i - \Lambda \eta_i + \Ad_{g_i^{-1}} p_\g \Ad_{g_i} \xi_i,
\;\;\;
z_i + \zeta_i = \Lambda(\eta_i - \xi_i) + \eta_i - \xi_i.
\end{equation}
Thus by our description of $T_{(q_1, q_2, q_3)} \G_P$,
to show \eqref{eq-alpha-triple-ok}, it suffices to show
\begin{align}
\lb{eq-22}
&\bar{x}_2:=x_2 - \Ad_{g_2^{-1}} p_\g 
\Ad_{\scriptscriptstyle{g_1^{l_1} l_1^{-1}}}(x_1 + z_1 + \zeta_1) \in \h,\\
\lb{eq-23}
&x_1 + \Ad_{h_3} x_3 \in \h\\
\lb{eq-33}
&\Ad_{h_3} (x_3 + z_3 + \zeta_3) +x_1 +z_1 + \zeta_1 +
\Ad_{l_1h}(\bar{x}_2 + z_2 + \zeta_2) \in \Ad_{l_1hl_2} \h.
\end{align} 
By \eqref{eq-i}, 
\begin{equation}\lb{eq-44}
x_i + z_i + \zeta_i = y_i + \eta_i -\Ad_{g_i^{-1}}
\Ad_{g_i^{-1}}^*\xi_i, \hs i = 1, 2, 3.
\end{equation}
We first prove  \eqref{eq-22}. Since $\Ad_{l_1^{-1}} (y_1 + \eta_1) \in \g$ and
$g_1^{l_1}l_1^{-1}g_1^{-1} \in G^*$ by \leref{le-in-Gs}, 
\[
\bar{x}_2 = y_2 -\Lambda \eta_2 + \Lambda \xi_2 
+\Ad_{g_2^{-1}}p_\g \Ad_{g_2} \xi_2 - \Ad_{h^{-1} l_1^{-1}}(y_1 + \eta_1).
\]
Note from \eqref{eq-2} that 
\[
\Ad_{h^{-1} l_1^{-1}}(y_1 + \eta_1) + \xi_2 -\Ad_{g_2^{-1}}p_\g \Ad_{g_2} \xi_2 +
\Ad_{h^{-1}l_1^{-1}h_3}(y_3 + \eta_3) \in \l,
\]
so by \eqref{eq-4},
$\Ad_{g_2^{-1}}p_\g \Ad_{g_2} \xi_2 - 
\Ad_{h^{-1} l_1^{-1}}(y_1 + \eta_1) -\xi_2 + y_2 + \eta_2 \in \l.$
It follows  from \eqref{eq-l-Lambda} that $\bar{x}_2 \in \h$, so  
\eqref{eq-22} holds.
To prove \eqref{eq-23}, note that \eqref{eq-1} implies that
$p_{\g^*} \Ad_{l_1 (g_1^{l_1})^{-1}} \Ad_{g_2^{-1}}^* \xi_2 = 
-\xi_1 - \Ad_{h_3^{-1}}^*\xi_3$, so
by \eqref{eq-2} and \eqref{eq-l-Lambda},
\[
y_1 + p_\g \Ad_{l_1 (g_1^{l_1})^{-1}} \Ad_{g_2^{-1}}^* \xi_2 + \Ad_{h_3} y_3 
+ p_\g \Ad_{h_3} \eta_3 +\Lambda (\xi_1-\eta_1) + \Lambda \Ad_{h_3^{-1}}^* (\xi_3 
-\eta_3) \in \h.
\]
Thus by \eqref{eq-1}, \eqref{eq-23} is equivalent to 
\[
\Ad_{h_3} \Lambda(\xi_3 - \eta_3) + p_\g \Ad_{h_3} (\xi_3 - \eta_3) -
\Lambda \Ad_{h_3^{_1}}^* (\xi_3 - \eta_3) \in \h
\]
which holds because of \leref{le-xi-h}. This proves \eqref{eq-23}. 
It remains to prove \eqref{eq-33}. Using \leref{le-in-Gs} and
\eqref{eq-44}, one sees that the left hand side of \eqref{eq-33} is equal to
\[
\Ad_{h_3} (y_3 + \eta_3) + \Ad_{l_1h} (y_2 + \eta_2) -
\Ad_{g_1^{-1}} p_{\g^*} \Ad_{g_1} (\xi_1 + \Ad_{h_3} \xi_3) - \Ad_{l_1 (g_1^{l_1})^{-1}}
\Ad_{g_2^{-1}}^* \xi_2.
\]
By \eqref{eq-1}, 
$\Ad_{g_1^{-1}} p_{\g^*} \Ad_{g_1} (\xi_1 + \Ad_{h_3} \xi_3) + \Ad_{l_1 (g_1^{l_1})^{-1}}
\Ad_{g_2^{-1}}^* \xi_2 = 0$. Moreover, since $y_2 + \eta_2 \in \Ad_{l_2} \g$ and
$y_3 + \eta_3 \in \Ad_{h_{3}^{_1}l_1hl_2} \g$,  
\[
y_2 + \eta_2 + \Ad_{h^{-1}l_1^{-1} h_3} (y_3 + \eta_3) \in \Ad_{l_2} \g \cap 
l = \Ad_{l_2} \h
\]
by \eqref{eq-4}. Thus the left hand side of \eqref{eq-33} is in
$\Ad_{l_1hl_2} \h$, so \eqref{eq-33} holds.

This finishes the proof of \thref{th-poi-grpoid}.

\end{document}